\pdfoutput=1
\makeatletter
\IfFileExists{./numapde-latex/numapde-packages.sty}{\providecommand*{\input@path}{}\edef\input@path{{./numapde-latex/}\input@path}}{}
\makeatother

\documentclass{numapde-preprint}

\usepackage{numapde-semantic}
\usepackage{numapde-style}
\usepackage{numapde-local}

\IfFileExists{./numapde-bibliography/numapde.bib}{\addbibresource{./numapde-bibliography/numapde.bib}}{\addbibresource{numapde.bib}}

\hypersetup{
	pdftitle={Continuous Galerkin Schemes for Semi-Explicit Differential-Algebraic Equations},
	pdfauthor={Robert Altmann, Roland Herzog},
	pdfkeywords={continuous Galerkin, Petrov-Galerkin schemes, semi-explicit DAE, variational consistency}
}

\title{Continuous Galerkin Schemes for Semi-Explicit Differential-Algebraic Equations\thanks{R.~Altmann acknowledges funding by the Deutsche Forschungsgemeinschaft (DFG, German Research Foundation) -- Project-ID~446856041.}}
\subtitle{}
\shorttitle{cG Schemes for Semi-Explicit DAEs}

\author{Robert Altmann\thanks{Department of Mathematics, University of Augsburg, Universitätsstraße~14, 86159~Augsburg, Germany (\email{robert.altmann@math.uni-augsburg.de}, \url{https://www.uni-augsburg.de/de/fakultaet/mntf/math/prof/numa/team/robert-altmann/}, \orcid{0000-0002-4161-6704}).}
\and
Roland Herzog\thanks{Interdisciplinary Center for Scientific Computing, Heidelberg University, 69120 Heidelberg, Germany (\email{roland.herzog@iwr.uni-heidelberg.de}, \url{https://www.tu-chemnitz.de/mathematik/part_dgl/people/herzog}, \orcid{0000-0003-2164-6575}).}}
\shortauthor{R. Altmann and R. Herzog}

\dedication{}

\IfFileExists{numapde-uncolorizeHyperrefIfChanges.sty}{\RequirePackage{numapde-uncolorizeHyperrefIfChanges}}{}

\begin{document}
\maketitle

% Insert abstract
\begin{abstract}
This paper studies a new class of integration schemes for the numerical solution of semi-explicit differential-algebraic equations of differentiation index~2 in Hessenberg form.
Our schemes provide the flexibility to choose different discretizations in the differential and algebraic equations. 
At the same time, they are designed to have a property called variational consistency, \ie, the choice of the discretization of the constraint determines the discretization of the Lagrange multiplier. 
For the case of linear constraints, we prove convergence of order~$r+1$ both for the state and the multiplier if piecewise polynomials of order~$r$ are used. 
These results are also verified numerically.\end{abstract}

% Insert keywords

% Insert Mathematics Subject Classification (MSC2010)
\begin{AMS}
\href{https://mathscinet.ams.org/msc/msc2010.html?t=65L80}{65L80}, \href{https://mathscinet.ams.org/msc/msc2010.html?t=65L60}{65L60}
\end{AMS}

% Insert document body
%=============================================================================
%=========  Sect 1: Introduction
%=============================================================================
\section{Introduction}
In this paper, we consider a novel class of integration schemes for the numerical solution of semi-explicit differential-algebraic equations (DAEs) of differentiation index~2 in Hessenberg form, 
\begin{subequations}
	\label{eq:DAE}
	\begin{align}
		\dot x 
		& 
		= 
		f(x,t) - g_x(x,t)^\transp \lambda
		, 
		\label{eq:DAE:a}
		\\
		0 
		& 
		= 
		g(x,t) 
		\label{eq:DAE:b}
	\end{align}
\end{subequations}
with an initial condition $x(0) = x_0$. 
Systems of type \eqref{eq:DAE} arise, for instance, as constrained gradient flow problems.
Indeed, consider the minimization of a possibly time-dependent energy $E(x,t)$ \wrt~$x$, subject to the constraint $g(x,t) = 0$.
The Karush--Kuhn--Tucker conditions associated with this problem read
\begin{equation*}
	0 
	= 
	\nabla_x E(x,t) + g_x(x,t)^\transp \lambda
	\quad \text{and} \quad 
	g(x,t) 
	= 
	0
	.
\end{equation*}
Hence, the corresponding Euclidean, constrained gradient flow is governed by a system of the form \eqref{eq:DAE} with $f = - \nabla_x E$. 
Systems of the form~\eqref{eq:DAE} also arise in fluid mechanics, where the constraint represents the incompressibility condition, cf.~the semi-discrete Oseen or (Navier-)Stokes equations; see \cite{Temam:1984:1}. 
Generally, time-dependent partial differential equations of first order (in time), which meet an additional constraint, lead to \eqref{eq:DAE} under spatial discretization~\cite{Altmann:2015:1}. 

The most popular numerical methods for systems of the form~\eqref{eq:DAE} are certainly (stiffly accurate) Runge-Kutta methods~\cite{HairerLubichRoche:1989:1,AltmannZimmer:2018:1}, and backward differentiation formulas, cf.~\cite[Ch.~3]{BrenanCampbellPetzold:1996:1}. 
For further approaches and details we refer to the surveys in~\cite[Ch.~VII]{HairerWanner:1996:1}, \cite[Ch.~10]{AscherPetzold:1998:1}, and~\cite[Ch.~5]{KunkelMehrmann:2006:1}. 
In the case of linear constraints, one may also consider splitting schemes~\cite{AltmannOstermann:2017:1}, discontinuous Galerkin schemes~\cite{VoulisReusken:2019:1,Voulis:2019:1}, or exponential integrators~\cite{AltmannZimmer:2020:1}. 
All three of these classes of methods were analyzed in the infinite-dimensional operator setting. 

In this paper, we develop a class of Petrov-Galerkin schemes, also known as continuous Galerkin (cG) methods, which can be directly applied to the DAE~\eqref{eq:DAE} without the need of an index reduction or regularization.  
To the best of our knowledge, this class of integrators has previously been applied and analyzed exclusively for ordinary differential equations (ODEs) and operator differential equations, \ie, in the absence of the terms involving the constraint function $g$; see for instance \cite{EstepFrench:1994:1,AkrivisMakridakis:2004:1,Wihler:2005:1,Schieweck:2010:1,MeidnerVexler:2011:1}. 
For an ODE~$\dot x = f(x,t)$, cG methods are based on the variational formulation
\begin{equation}
	\label{eq:variational_formulation_ODE}
	\int_0^T [\dot x(t) - f(x(t),t)]^\transp w(t) \dt 
	= 
	0
	.
\end{equation}
They are obtained when $x$ is discretized as a vector-valued, globally continuous, piecewise polynomial function of degree~$r$ on a partitioning of the solution interval~$[0,T]$, while the test functions~$w$ are vector-valued, piecewise polynomials of degree~$r-1$ but discontinuous on the same partitioning. 
A specific scheme is then essentially determined by the quadrature formula applied to handle the term involving the right hand side~$f$, which is generally nonlinear. The choice of a basis in the solution and test spaces is also relevant from a practical point of view but it does not affect the scheme per se. 
A number of well-known one-step methods for ODEs can be generated in this way, including the trapezoidal rule  as the lowest-order member, which is obtained when~$r = 1$ and the trapezoidal rule is used for the quadrature of $f(x(t),t)^\transp w(t)$ in~\eqref{eq:variational_formulation_ODE}. 
This scheme is also known as the Crank-Nicolson method, in particular when applied to time-dependent partial differential equations.

In this paper, we extend cG schemes to semi-explicit DAEs \eqref{eq:DAE}, \ie, we need to address the discretization of the constraint $g(x,t) = 0$ as well as the additional forcing term $g_x(x,t)^\transp \lambda$.
Our guideline is a principle which we term \emph{variational consistency}.
To explain it, consider once more the case when $f$ is a (possibly time-dependent) gradient field, \ie, $f(x,t) = - \nabla_x E(x,t)$ holds.
In this case, \eqref{eq:DAE} can be viewed as a continuous gradient flow subject to the constraint $g(x,t) = 0$ and $\lambda$ is the corresponding Lagrange multiplier.
When one discretizes \eqref{eq:DAE} in time, the choices how to treat the terms $g$ and $g_x^\transp \lambda$ appear to be independent. 
This philosophy is known as \emph{first optimize, then discretize}, compare for instance the discussion in the context of optimal control for partial differential equations in \cite[Ch.~3.2]{HinzePinnauUlbrichUlbrich:2009:1}.
By contrast, we honor the variational consistency of our family of schemes by following the \emph{first discretize, then optimize} approach. 
As a consequence, two choices completely determine the scheme. 
We only need to choose a quadrature formula which determines the discretization of the ODE part described by~$f$, and to decide in which way to enforce the constraint~$g = 0$. 
The discretization of the remaining term $g_x^\transp \lambda$ will then follow automatically by an evaluation of the optimality conditions associated with the discretized problem. 
The resulting numerical schemes and a detailed discussion on the particular choice of the polynomial basis is part of~\Cref{sec:development}. 
We would like to emphasize that the principle of variational consistency applies also when $f$ is not a gradient field. 
Moreover, although this is not fully explored in the present paper, our methods are naturally fully adaptive in terms of the local time step size and polynomial degree (\emph{$hp$-adaptive}).

Due to the different treatment of the differential and algebraic equations in~\eqref{eq:DAE}, we obtain a new class of high-order DAE integrators. 
Note that this flexibility calls for a concept such as variational consistency in the first place. 
A second difference compared to standard schemes such as Runge-Kutta methods is the distinction of the discrete ansatz spaces for the variables~$x$ and~$\lambda$. 
Note that this is also the case in collocation methods, where the Lagrange multiplier is typically approximated by a piecewise polynomial of one degree lower than~$x$, cf.~\cite{HankeMaerzTischendorfWeinmuellerWurm:2017:1}. 
Here, however, we approximate~$\lambda$ by a linear combination of functionals. 

For the case of linear constraints (\wrt~$x$), we show in~\Cref{sec:errAnalysis} that the introduced class of cG schemes indeed leads to high-order integrators. 
More precisely, the rates known for cG schemes from the ODE case are recovered and they also extend to the Lagrange multiplier. 
Again the analysis is simplified significantly by the concept of variational consistency, since the saddle point structure is maintained. 
Finally, the convergence properties are verified by a number of numerical experiments in~\Cref{sec:numerics}. 
The examples include a simple linear circuit problem and a more challenging coupled quasilinear heat equation with a thermal resistance condition between two subdomains.
We also include results for the pendulum problem, which has a similar structure than \eqref{eq:DAE} but contains second-order derivatives.
It is therefore not covered by the present analysis but still our scheme is applicable and convergent.

The code for this paper is available on GitHub\footnote{https://github.com/rolandherzog/cg-schemes-for-daes} and can be cited as \cite{AltmannHerzog:2021:1}.

%=============================================================================
%=========  Sect 2: cG Schemes
%=============================================================================
\section{Development of Continuous Galerkin Schemes}
\label{sec:development}
In this section, we develop the class of variationally consistent continuous Galerkin schemes for \eqref{eq:DAE}.
We seek an approximation of the solution $(x,\lambda) \colon [0,T] \to \R^n \times \R^m$ to~\eqref{eq:DAE} on a given time horizon $[0,T]$, subject to the initial conditions $x(0) = x_0 \in \R^n$. 
We call the initial data consistent if it is in line with the constraint~\eqref{eq:DAE:b}, \ie, $g(x_0,0) = 0$. 
Such a condition is necessary for the existence of continuous solutions but not for the construction of the proposed cG schemes.  
The right-hand side and the constraint function are given as $f \colon \R^n \times [0,T] \to \R^n$ and $g \colon \R^n \times [0,T] \to \R^m$.
We make the following standing assumption.
\begin{assumption}
	\label{assumption:standing_assumption}
	Both $f$ and $g$ are assumed to be sufficiently smooth and the derivative of~$g$ \wrt~$x$ should be of full row rank. 
	In particular, we assume that~\eqref{eq:DAE} has a unique solution with initial condition~$x(0)=x_0$.
\end{assumption}

Let us denote by $\cP_r(I)$ and $\cP_r(I;\R^n)$ the space of scalar and vector-valued polynomials, respectively, of maximal degree $r \ge 0$ on some interval $I \subset \R$.
As already mentioned in the introduction, cG schemes employ a continuous solution space
\begin{equation*}
	V 
	\coloneqq 
	\setDef[auto]{x \in C([0,T];\R^n)}{\restr{x}{I_\ell} \in \cP_r(I_\ell; \R^n)}
	.
\end{equation*}
Here $I_1, \ldots, I_N$ denotes a partitioning of $(0,T]$ into intervals $I_\ell \coloneqq (T_{\ell-1},T_\ell]$ of lengths $\Delta_\ell \coloneqq T_\ell - T_{\ell-1}$.
Since the test space 
\begin{equation*}
	W 
	\coloneqq 
	\setDef[auto]{w \colon [0,T] \to \R^n}{\restr{x}{I_\ell} \in \cP_{r-1}(I_\ell; \R^n)}
\end{equation*}
consists of discontinuous functions, variational formulations such as \eqref{eq:variational_formulation_ODE} decouple into individual time intervals, allowing \eqref{eq:variational_formulation_ODE} to be written as a time-stepping method.
For notational convenience, we therefore concentrate on a single interval $I=[0,\Delta]$ of length~$\Delta$ and drop its index.
%
%=============================================================================
\subsection{Discretization of the ODE Part}
\label{sec:development:ODE}
Without loss of generality, we choose pairwise distinct Lagrange points $t_0, \ldots, t_r \in I$ with $t_0 = 0$. 
These points are not necessarily equispaced. 
We consider the Lagrange interpolation polynomials $\varphi_0, \ldots, \varphi_r$ as our basis in $\cP_r(I)$, the ansatz space for the state. 
In other words, $\varphi_j(t_k) = \delta_{j,k}$ holds for $0 \le j,k \le r$.
Then the approximation of the state $x$ on the interval $I$ can be written as 
\begin{equation*}
	\x(t) 
	= 
	\sum_{j=0}^r x_j \, \varphi_j(t)
	, 
	\quad 
	t \in I
\end{equation*}
with unknown coefficients $x_j \in \R^n$. 
Due to the choice of a Lagrangian basis, $x_j = X(t_j)$ holds. 
We denote the basis for the test space $\cP_{r-1}(I)$ by $\psi_1, \ldots, \psi_r$. 
A convenient choice of this basis will become evident later.

We begin with the discussion of the discretization of the ODE part of \eqref{eq:DAE}, \ie, $\dot x - f(x,t) = 0$, on the interval $I$.
Replacing $f$ by its nodal interpolation in the same Lagrange points as above, \ie, 
\begin{equation*}
	f(x(t),t)
	\approx
	\sum_{j=0}^r f(x_j,t_j) \, \varphi_j(t)
	,
\end{equation*}
and utilizing the bases for the solution and test spaces, we can approximate the variational formulation \eqref{eq:variational_formulation_ODE} on $I$ by 
\begin{equation}
	\label{eq:cG_scheme_ODE_part}
	\sum_{j=0}^r x_j \int_I \dot \varphi_j(t) \, \psi_i(t) \dt
	-
	\sum_{j=0}^r f(x_j,t_j) \int_I \varphi_j(t) \, \psi_i(t) \dt
	=
	0
	,
\end{equation}
where $i = 1, \ldots, r$.
Notice that these are $r\, n$ scalar, generally nonlinear equations for the same number of unknowns $x_1, \ldots, x_r \in \R^n$, since the initial value $x_0 = x(t_0) = x(0) \in \R^n$ on the current interval coincides with the terminal value of the previous interval, or is specified through the initial conditions.

The material up to here is standard in the development of cG methods for ODEs; see for instance~\cite[Ch.~9.2.1]{ErikssonEstepHansboJohnson:1996:1} or~\cite{EstepFrench:1994:1,Schieweck:2010:1}.
Let us now focus on the novel treatment of the constraint term $g(x,t) = 0$ and the induced, variationally consistent discretization of $g_x(x,t)^\transp \lambda$. 
%
%=============================================================================
\subsection{Discretization of the Constraint and Variational Consistency}
\label{subsection:variational_consistency}
Clearly, it is in general impossible to satisfy the constraint $g(x(t),t)$ everywhere on the time interval~$I$.
We therefore enforce this constraint in~$r$ pairwise distinct time points, namely $s_1, \ldots, s_r \in I$. 
Since $x(0)$ is already fixed, it is reasonable to assume that $0 = t_0 \neq s_k$ for any $k = 1, \ldots, r$. 
The corresponding functionals of point evaluation (Dirac measures) are denoted by $\xi_1, \ldots, \xi_r$.
That is, by definition, we have
\begin{equation*}
	\dual{\xi_k}{g(x(\cdot), \cdot)}
	\coloneqq
	g(x(s_k), s_k)
	,
	\quad 
	k = 1, \ldots, r
	.
\end{equation*}
Hence, we do not consider polynomials as test space for the constraint but functionals spanned by $\xi_1, \ldots, \xi_r$. 
More precisely, we consider the test functionals 
\begin{equation*}
	\Lambda 
	= 
	\sum_{k=1}^r \lambda_k \, \xi_k
\end{equation*}
with coefficients $\lambda_k\in\R^m$. For the discrete setting this then leads to 
\begin{equation}
	\label{eq:understanding_of_constraint}
	0 
	= 
	\dual{\Lambda}{g(x(\cdot), \cdot)}
	= 
	\sum_{k=1}^r \lambda_k^\transp \dual{\xi_k}{g(x(\cdot), \cdot)}
	= 
	\sum_{k=1}^r \lambda_k^\transp g(x(s_k), s_k)
	\quad \text{for all } \lambda_k \in \R^m
	,
\end{equation}
which is equivalent to $0=g(x(s_k), s_k)$ for $k=1, \ldots, r$. 

It remains to discuss the discretization of the term $g_x(x,t)^\transp \lambda$. Following the mentioned philosophy known as \emph{first discretize, then optimize}, this discretization follows automatically by evaluating the optimality conditions associated with the discretized problem. 
Thus, we consider the directional derivative of~\eqref{eq:understanding_of_constraint} \wrt~$x$ in the direction of a basis function $\psi_i$ of the local test space $\cP_{r-1}(I)$.
Applying the chain rule, this amounts to
\begin{equation*}
	\restr[auto]{\ddx \paren[Bigg](){ \sum_{k=1}^r \lambda_k^\transp g(x(s_k), s_k) }}{\psi_i}
	= 
	\sum_{k=1}^r g_x(x(s_k), s_k)^\transp \lambda_k \, \psi_i(s_k)
	.
\end{equation*}

Up to now, we have only fixed the ansatz functions $\varphi_0, \ldots, \varphi_r$ as the Lagrange basis functions \wrt the points $0 = t_0, \ldots, t_r$ and we decided to use the same points for the quadrature rule applied to~$f$. 
However, we are still free to specify the basis functions~$\psi_i$ as well as the evaluation points $s_k$ defining the functionals~$\xi_k$. 
In this general setting, the Galerkin scheme has the following form: 
on the interval~$I$, find 
\begin{equation*}
	\x(t) 
	= 
	\sum_{j=0}^r x_j \, \varphi_j(t)
	, 
	\qquad 
	\Lambda 
	= 
	\sum_{k=1}^r \lambda_k \, \xi_k
	,
\end{equation*}
where the unknown coefficients $x_1, \ldots, x_r \in \R^n$ and $\lambda_1, \ldots, \lambda_r \in \R^m$ are the solution of the nonlinear system on the interval $I$:
\begin{subequations}
	\label{eq:cG_time_step}
	\begin{align}
		\sum_{j=0}^r \alert{x_j} \, D_{ij}
		-
		\sum_{j=0}^r f(\alert{x_j},t_j) \, M_{ij}
		+
		\sum_{k=1}^r g_x\paren[Big](){\sum_{j=0}^r\alert{x_j} \, \varphi_j(s_k), s_k}^\transp \alert{\lambda_k}\, \psi_i(s_k)
		&
		=
		0
		\quad
		\text{for } i = 1,\ldots, r,
		\label{eq:cG_time_step_1}
		\\
		g\paren[Big](){\sum_{j=0}^r\alert{x_j} \, \varphi_j(s_k),s_k} 
		&
		=
		0
		\quad
		\text{for } k = 1,\ldots, r.
		\label{eq:cG_time_step_2}
	\end{align}
\end{subequations}
Here and throughout, we make use of the definitions
\begin{subequations}
	\label{eq:defDM}
	\begin{align}
		(D_{ij}) 
		\coloneqq
		\left( \int_I \dot \varphi_j(t) \, \psi_i(t) \dt \right)
		&
		\in \R^{r \times (r+1)}, 
		\\
		(M_{ij}) 
		\coloneqq
		\left( \int_I \varphi_j(t) \, \psi_i(t) \dt \right)
		&
		\in \R^{r \times (r+1)}.
	\end{align}
\end{subequations}
The partition of unity property (which is satisfied for Lagrange polynomials) yields the following result on the row and column sums of $D$. 

\begin{lemma}
	\label{lemma:row_column_sums_of_Dbar}
	Suppose that the basis functions $\psi_1, \ldots, \psi_r$ satisfy the property $\sum_{i=1}^r \psi_i \equiv 1$. 
	Further assume that~$t_r = \Delta$. 
	Then, for any $r \ge 1$, the matrix $D$ has the property 
	\begin{equation*}
		\sum_{i=1}^r D_{ij} = 
		\begin{cases}
			-1 & \text{if } j = 0 
			\\
			0 & \text{if } 1 \le j \le r-1 
			\\
			1 & \text{if } j = r
		\end{cases}
		\quad \text{and} \quad
		\sum_{j=0}^r \, D_{ij} = 0
		\quad \text{for all } i = 1, \ldots, r.
	\end{equation*}
\end{lemma}
\begin{proof}
	By the definition of $D$, we have 
	\begin{equation*}
		\sum_{i=1}^r D_{ij} 
		=
		\sum_{i=1}^r \int_0^\Delta \dot \varphi_j \, \psi_i \, \dt
		=
		\varphi_j(\Delta) - \varphi_j(0),
	\end{equation*}
	since $\sum_{i=1}^r \psi_i \equiv 1$.
	This shows the first claim, since $t_0 = 0$ and $t_r = \Delta$. 
	Further we obtain 
	\begin{equation*}
		\sum_{j=0}^r \, D_{ij} 
		=
		\sum_{j=0}^r \int_0^\Delta \dot \varphi_j \, \psi_i \, \dt
		=
		0,
	\end{equation*}
	since $\sum_{j=0}^r \varphi_j \equiv 1$ and thus $\sum_{j=0}^r \dot \varphi_j \equiv 0$ holds.
\end{proof}

In the remainder of this section, we will use the given freedom in the choice of 
the basis $\{\psi_1, \ldots, \psi_r\}$ of test functions for \eqref{eq:cG_time_step} in $\cP_{r-1}(I)$ as well as in the points $s_1, \ldots s_r$, where the constraint is enforced, in order to optimize the structure of the resulting systems which have to be solved in every time step. 
%
%=============================================================================
\subsection{Optimizing the Structure}
\label{sec:development:structure}
System \eqref{eq:cG_time_step} is generally nonlinear unless both $f$ and $g$ are affine functions \wrt $x$. 
It is thus natural to solve it by Newton's method. 
Our choice of the basis functions $\psi_i$ and time points $s_k$ will be driven by the structure of the arising linear systems. Thus, we need to consider the derivative of \eqref{eq:cG_time_step} \wrt~the coefficients~$x_j$ and~$\lambda_k$. 

% (2,1)-block
We start with the derivative of the constraint equation~\eqref{eq:cG_time_step_2} \wrt~$x_j$. Since $x_0$ is not an unknown, we only need to consider $j = 1, \ldots, r$. 
This results in a block matrix~$\kronG$ where the entry $\kronG_{k, j}$ is a matrix of dimension $m\times n$ and  defined through
\begin{equation*}
	\kronG_{k, j}
	= \ddxj g\paren[Big](){\sum_{j=0}^r x_j \, \varphi_j(s_k),s_k}
	= g_x(X(s_k),s_k)\, \varphi_j(s_k) \in \R^{m\times n}. 
\end{equation*} 
% (1,2)-block
On the other hand, the derivative of equation~\eqref{eq:cG_time_step_1} \wrt~$\lambda_k$ yields a block matrix~$\hat \kronG$. 
Here, the entry $\hat \kronG_{i,\ell}$ is of dimension~$n \times m$ and equals
\begin{equation*}
	\hat \kronG_{i,\ell}
	= 
	\frac{\partial}{\partial \lambda_\ell} \sum_{k=1}^r g_x(X(s_k), s_k)^\transp \lambda_k\, \psi_i(s_k)
	= 
	g_x(X(s_\ell), s_\ell)^\transp\, \psi_i(s_\ell) \in \R^{n\times m}
	.
\end{equation*} 
These block matrices $\kronG$ and $\hat \kronG$ define the off-diagonal blocks in the Newton matrix. 
The choice $s_k = t_k$ for $k = 1, \ldots, r$, \ie, enforcing the constraint exactly in the Lagrange points $t_1, \ldots, t_r$, renders the matrix~$\kronG$ block-diagonal with blocks $\kronG_{k,k} = g_x(x_k,t_k)$. 
Further, we may fix the basis functions $\psi_1, \ldots, \psi_r$ such that
\begin{equation*}
	\psi_i(s_k)
	= \psi_i(t_k)
	= \delta_{i,k}.
\end{equation*} 
This means that the $\psi_i$ are the Lagrange basis of $\cP_{r-1}(I)$ in the points $t_1, \ldots, t_r$. 
With this particular choice, also the matrix $\hat \kronG$ is block-diagonal and we have $\hat \kronG = \kronG^\transp$. 

With the choices of this subsection, system~\eqref{eq:cG_time_step} for the unknowns $\alert{x_1}, \ldots, \alert{x_r}$ and $\alert{\lambda_1}, \ldots, \alert{\lambda_r}$ simplifies to
\begin{subequations}
	\label{eq:cG_time_step_simplified}
	\begin{align}
		\sum_{j=0}^r \alert{x_j} \, D_{ij}
		-
		\sum_{j=0}^r f(\alert{x_j},t_j) \, M_{ij}
		+
		g_x(\alert{x_i}, t_i)^\transp \alert{\lambda_i}
		&
		=
		0
		\quad
		\text{for } i = 1,\ldots, r,
		\label{eq:cG_time_step_simplified_1}
		\\
		g(\alert{x_k},t_k)
		& 
		=
		0
		\quad
		\text{for } k = 1,\ldots, r.
		\label{eq:cG_time_step_simplified_2}
	\end{align}
\end{subequations}
This formulation allows us to interpret the Lagrange multiplier as a weighted sum of point forces in the Lagrange points $t_1, \ldots, t_r$, in which the constraint~$g = 0$ is enforced. 
Recall that the resulting approximation $\Lambda$ was defined in terms of the basis $\xi_k$, which consists of functionals.

We emphasize that~\eqref{eq:cG_time_step_simplified} still contains the freedom of how to choose the Lagrange points $t_1, \ldots, t_r$ for enforcement of the constraint. 
This choice also determines the basis functions~$\varphi_j$ and~$\psi_i$ and therefore the matrices $D$ and $M$. The following result is independent of this choice. 
\begin{theorem}
	\label{thm:solvability}	
	Assume that $t_0, \ldots, t_r$ are pairwise distinct Lagrange points and that $\Delta > 0$ is a sufficiently small step size. 
	Further assume that the constraint is linear \wrt~$x$ in the sense that~$g$ has the particular form~$g(x,t) = \tilde Bx - \tilde g(t)$. 
	Then system~\eqref{eq:cG_time_step_simplified} has at least locally a unique solution, \ie, the cG scheme is well-defined. 
\end{theorem}
We postpone the proof to the end of this section and first discuss some examples. 
%
%
%=============================================================================
\subsection{Examples}
\label{sec:development:examples}
Given the assumptions from the previous subsection, we consider the resulting schemes for the lowest-degree cases~$r=1$ and $r=2$, \ie, for piecewise linear and quadratic approximations of $x$. 
In each case, we consider a single time step on an interval~$I$ of length~$\Delta$.
For demonstration purposes, we choose $t_0, \ldots, t_r$ equispaced so that $\varphi_0, \ldots, \varphi_r$ becomes the standard Lagrange basis. 

\subsubsection{Case \texorpdfstring{$r=1$}{r=1}}
The variable $x$ is approximated by an affine function, where the value~$x_0$ at the left end of the interval is fixed from the previous step or the initial condition. 
Thus, the unknowns only consist of $x_1 \in \R^n$ and~$\lambda_1 \in \R^m$, which define 
\begin{equation*}
	X(t) = x_0 \varphi_1(t) + x_1 \varphi_1(t)
	, 
	\qquad 
	\Lambda = \lambda_1  \xi_1.  
\end{equation*}
In this case, we have~$D = [-1,\, 1]$ and $M = \tfrac \Delta2\, [1,\, 1]$ and the cG scheme~\eqref{eq:cG_time_step_simplified} turns into
\begin{subequations}
	\label{eq:cG_r0}
	\begin{align}
		\alert{x_1} - x_0   
		- \tfrac{\Delta}{2} \, \paren[big][]{f(x_0,t_0) + f(\alert{x_1},t_1)}
		+ g_x(\alert{x_1},t_1)^\transp \alert{\lambda_1} 
		&
		=
		0
		, 
		\label{eq:cG_r0_1} 
		\\
		g(\alert{x_1},t_1)
		&
		= 
		0
		.
		\label{eq:cG_r0_2}
	\end{align}
\end{subequations}
Thus, we obtain the trapezoidal rule for the ODE part.
However, the constraint is included implicitly only at the end of the interval, at $t_1 = \Delta$.
The principle of variational consistency then determines how the term $g_x(x,t)^\transp \lambda$ is discretized.

A Newton step for \eqref{eq:cG_r0} is governed by the matrix
\begin{equation*}
	\begin{bmatrix}
		\id_n - \frac{\Delta}{2} f_x(x_1,t_1) 
		& 
		g_x(x_1,t_1)^\transp 
		\\
		g_x(x_1,t_1) 
		& 
		0
	\end{bmatrix}
	,
\end{equation*}
where $\id_n$ denotes the $n \times n$-identity matrix.

\begin{remark}
	The	numerical scheme~\eqref{eq:cG_r0} cannot be written as a Runge-Kutta scheme. 
	To see this, we consider a direct application of the Butcher tableau 
	\begin{equation*}
		\begin{array}{r|cc}
			0 &     &     \\
			1 & 1/2 & 1/2 \\
			\hline
			& 1/2 & 1/2
		\end{array} 
	\end{equation*}
	to~\eqref{eq:DAE}.
	This corresponds to the trapezoidal rule for ODEs and it results in update rule similar to~\eqref{eq:cG_r0} but with a second constraint $g(x_0,t_0) = 0$ and a corresponding second stage of the Lagrange multiplier. 
	Note, however, that this second constraint is a condition for the (given) value $x_0$, which is naturally satisfied for a consistent value $x_0$ or unsolvable otherwise. 
	The different nature of the cG scheme developed in this paper comes from the fact that we consider different discretizations for~$f$ and~$g$. 
\end{remark}
\subsubsection{Case \texorpdfstring{$r=2$}{r=2}}
The solution component $x$ is now approximated by a quadratic function, and the unknowns are the coefficients~$x_1, x_2 \in \R^n$. 
The Lagrange multiplier is described by coefficients $\lambda_1, \lambda_2\in \R^m$. 
In the present case, the coefficients in~\eqref{eq:cG_time_step_simplified} are given by 
\begin{equation*} 
	D = \frac 13 \begin{bmatrix} -5 & 4 & 1 \\ 2 & -4 & 2 \end{bmatrix}, \qquad 
	M = \frac \Delta 6 \begin{bmatrix} 2 & 4 & 0 \\ -1 & 0 & 1 \end{bmatrix}. 
\end{equation*}
The resulting approximations are hence given by 
\begin{equation*}
	X(t) = x_0 \varphi_0(t) + x_1 \varphi_1(t) + x_2 \varphi_2(t), \qquad 
	\Lambda = \lambda_1 \xi_1 + \lambda_2 \xi_2  
\end{equation*}
with coefficients solving 
\begin{align*}
	- \frac{5}{3} x_0 + \frac{4}{3}\alert{x_1} + \frac{1}{3}\alert{x_2}
	- \frac{\Delta}{6} \, \big[ 2f(x_0, t_0) + 4f(\alert{x_1}, t_1) \big]
	+ g_x(\alert{x_1}, t_1)^\transp \alert{\lambda_1} 
	&
	= 
	0
	, 
	\\
	\frac{2}{3} x_0 - \frac{4}{3} \alert{x_1} + \frac{2}{3} \alert{x_2}
	- \frac{\Delta}{6} \, \big[ -f(x_0, t_0) + f(\alert{x_2}, t_2) \big] 
	+ g_x(\alert{x_2}, t_2)^\transp \alert{\lambda_2}
	&
	= 
	0
	, 
	\\
	g(\alert{x_1}, t_1) 
	= 
	0
	\quad \text{and} \quad
	g(\alert{x_2}, t_2) 
	&
	= 
	0
	. 
\end{align*}
We conclude this section with a discussion of the solvability of the cG scheme as well as the practical implementation of the method. 

%
%=============================================================================
\subsection{Newton Matrices and Solvability}
\label{sec:development:Newton}
We now address the numerical solution of the nonlinear problem \eqref{eq:cG_time_step_simplified} to be solved on each interval, which we denote as $I = [0,\Delta]$ for convenience.
Recall that $t_0, \ldots, t_r \in I$ are pairwise distinct Lagrange points on $I$, not necessarily equispaced. 
\begin{assumption}
	\label{assumption:left_and_right_Lagrange_points}
	In addition to $t_0 = 0$ we now also assume $t_r = \Delta$. 
\end{assumption}
Notice that, due to the continuity of the state approximation~$X$, \Cref{assumption:left_and_right_Lagrange_points} implies that the initial value for the subsequent time interval will be consistent.

Our method of choice for the solution of \eqref{eq:cG_time_step_simplified} is Newton's method.
In order to describe a Newton step, we introduce the following notation.
The unknowns located on the interval~$I$ are denoted as 
\begin{equation*}
	\xVec 
	\coloneqq 
	\begin{bmatrix}
		x_1 \\ \vdots \\ x_r
	\end{bmatrix}
	\in \R^{r n}
	\quad \text{and} \quad
	\lVec \coloneqq
	\begin{bmatrix}
		\lambda_1 \\ \vdots \\ \lambda_r
	\end{bmatrix}
	\in \R^{r m}
	.
\end{equation*}
We partition $D$ and $M$ into
\begin{equation}
	\label{eq:partitioning_D_and_M}
	D = 
	\begin{bmatrix}
		D_1 & \overline D
	\end{bmatrix}
	\quad \text{and} \quad
	M =
	\begin{bmatrix}
		M_1 & \overline M
	\end{bmatrix}
	,
\end{equation}
where $D_1$, $M_1$ are the first columns and $\overline D$ and $\overline M$ are square of size $r \times r$.
Notice that $\overline D$ and $\overline M$ are the matrices pertaining to the unknowns $x_1, \ldots, x_r$ in \eqref{eq:cG_time_step_simplified}.
\begin{remark}
	\label{remark:properties_of_Dbar}
	It can be observed that the matrix $\overline D$ is generally non-symmetric for $r \ge 2$.
	Interestingly, we observed numerically that the number of positive, negative, and zero eigenvalues of the \emph{symmetric part} $(\overline D + \overline D^\transp)/2$ are
	\begin{equation*}
		n^+ = \min\{r,2\},
		\qquad
		n^- = \max\{\min\{r-2,1\},0\},
		\qquad
		n^0 = \max\{r-3,0\}
	\end{equation*}
	for $r \ge 1$ when the points $t_0, \ldots, t_r \in I$ are equispaced; compare the examples given for $r \in \{1,2\}$ in \Cref{sec:development:examples}.
\end{remark}
As before, we denote by $\id_n$ the $n \times n$-identity matrix and define the Kronecker product matrices
\begin{equation*}
	\kronD 
	\coloneqq
	\overline D \otimes \id_n
	\quad \text{and} \quad
	\kronM 
	\coloneqq
	\overline M \otimes \id_n
\end{equation*}
of dimension $r n \times r n$.
The linearizations of the right hand side~$f$ and the constraint~$g$ give rise to the matrices
\begin{equation*}
	\begin{aligned}
		\kronF(\xVec)
		&
		\coloneqq
		\blkdiag\paren[big](){f_x(x_1,t_1), \ldots, f_x(x_r,t_r)}
		\in \R^{r n \times r n}
		,
		\\
		\kronG(\xVec)
		&
		\coloneqq
		\blkdiag\paren[big](){g_x(x_1,t_1), \ldots, g_x(x_r,t_r)}
		\in \R^{r m \times r n}
		.
	\end{aligned}
\end{equation*}

A Newton step for \eqref{eq:cG_time_step_simplified} at the iterate~$\xVec$ is governed by the block matrix
\begin{equation}
	\label{eq:Newton_matrix}
	\begin{bmatrix}
		\kronD - \kronM \kronF(\xVec) 
		& 
		\kronG(\xVec)^\transp
		\\
		\kronG(\xVec) 
		& 
		0
	\end{bmatrix}
	.
\end{equation}
Recall that the matrix $\kronM$ is proportional to the time step size~$\Delta$ while all other matrices are independent of~$\Delta$.

\begin{remark} \label{remark:we_cannot_get_symmetry_and_diagonality}
	Clearly, when $r = 1$ and $f(x,t) = \nabla_x E(x,t)$, then the (1,1) block and thus the entire matrix in \eqref{eq:Newton_matrix} is symmetric, which can be exploited in direct or iterative solvers.
	However, this is no longer the case when $r \ge 2$, since then $\overline D$ and thus $\kronD$ are generally non-symmetric, and $\overline M$ and thus $\kronM$ are generally non-diagonal and thus $\kronM$ does not commute with a general block-diagonal, symmetric~$\kronF(\xVec)$. 
	One might ask whether it is possible to choose interpolation points $t_1, \ldots, t_r$ which render $\overline D$ symmetric and $\overline M$ diagonal. 
	Unfortunately, a numerical experiment for $r = 2$ shows that this is not the case. 
	Indeed, a tabulation of $\overline D$ and $\overline M$ for numerous values of $0 = t_0 < t_1 < t_2 \le \Delta$ shows that $\overline M$ is diagonal for $(t_1,t_2) = (2/6,5/6)$ as well as for the equispaced points $(t_1,t_2) = (1/2,1)$, but $\overline D$ is never symmetric. 
\end{remark}

Our first result concerns the invertibility of the Newton matrix in \eqref{eq:Newton_matrix}.
The following lemma is an auxiliary step.
\begin{lemma}
	\label{lemma:Dbar_is_invertible}
	The matrix $\overline D$ is invertible.
\end{lemma}
\begin{proof}
	Recall that $\overline D$ has entries
	\begin{equation*}
		\overline D_{ij}
		=
		\int_0^\Delta \dot \varphi_j(t) \, \psi_i(t) \dt
		.
	\end{equation*}
	Due to the numbering convention, $i,j$ run from $1$ to $r$.
	We show that $\{\dot \varphi_1, \ldots, \dot \varphi_r\}$ are linearly independent and thus they span the space~$\cP_{r-1}(I)$.
	To see this, consider first
	\begin{equation*}
		0 
		=
		\sum_{j=1}^r \alpha_j \, \dot \varphi_j
		\quad
		\Rightarrow
		\quad
		0
		=
		\sum_{j=1}^r \int_0^{t_r} \alpha_j \, \dot \varphi_j \, \dt
		=
		\sum_{j=1}^r \alpha_j \, (\varphi_j(\Delta) - \varphi_j(0))
		=
		\alpha_r
		.
	\end{equation*}
	Repeating this argument and replacing the integration limit $t_r$ by $t_{r-1}$ etc.\ eventually shows $\alpha_r = \alpha_{r-1} = \ldots = \alpha_1 = 0$.
	Moreover, $\{\psi_1,\ldots,\psi_r\}$ constitutes another basis of $\cP_{r-1}(I)$.

	We infer that $\overline D$ contains the mutual $L^2$-inner products of two sets of basis vectors and is thus non-singular.
\end{proof}

Since $\kronD$ is defined as a Kronecker product of $\overline D$ and the identity matrix, \Cref{lemma:Dbar_is_invertible} directly implies the invertibility of $\kronD$. 
This, in turn, leads to the following result for the Newton matrix. 
\begin{proposition}
	\label{proposition:invertibility_Newton_matrix}
	Suppose that $\overline x \in \R^n$ is given, that $\Delta > 0$ and $\delta > 0$ are sufficiently small, and that $x_1, \ldots, x_r \in B_\delta(\overline x)$ holds.
	Then the Newton matrix in \eqref{eq:Newton_matrix} is nonsingular.
\end{proposition}
\begin{proof}
	Our proof utilizes a general criterion for non-symmetric saddle-point systems; see for instance \cite[Th.~3.1]{GanstererSchneidUeberhuber:2003:1}. 
	It is hence necessary and sufficient to show that $(i)$ $\kronG(\xVec)$ has full row rank and $(ii)$ $\kronZ^\transp \paren[big](){\kronD - \kronM \kronF(\xVec)} \kronZ$ is nonsingular, where~$\kronZ$ is a matrix with linearly independent columns spanning $\ker \kronG(\xVec)$.

	Condition~$(i)$ is satisfied by our standing~\Cref{assumption:standing_assumption}.
	In order to show $(ii)$, we are going to construct the matrix $\kronZ$ such that it is close to being block diagonal with repeated entries in order to exploit \Cref{lemma:Dbar_is_invertible}.
	To this end, let $\overline Z \in \R^{n \times (n-m)}$ be a matrix with linearly independent columns spanning $\ker g_x(\overline x,\Delta/2)$.
	By assumption, $x_j$ is close to $\overline x$ and $t_j$ is close to $\Delta / 2$ for $j = 1, \ldots, r$.
	Therefore, we can consider $g_x(x_j,t_j)$ to be small perturbations of $g_x(\overline x,\Delta/2)$.
	Consequently, we can find matrices $Z_j\in \R^{n \times (n-m)}$ spanning $\ker g_x(x_j,t_j)$, which are close to~$\overline Z$.
	More precisely, for any $\varepsilon > 0$ we can find $\Delta_0 > 0$ and $\delta_0 > 0$ such that for any $\Delta \in (0,\Delta_0)$ and any $\delta \in (0,\delta_0)$ and any collection of points $x_1, \ldots, x_r \in B_\delta(\overline x)$, there exist matrices $Z_j$ spanning $\ker g_x(x_j,t_j)$, $j = 1, \ldots, r$, such that $\norm{Z_j - \overline Z}_F < \varepsilon$ holds in the Frobenius norm, see~\cite{Dolezal:1964:1}. We can thus write $Z_j = \overline Z + \varepsilon \, E_j$ with $\norm{E_j}_F < 1$.

	We set $\kronZ \coloneqq \blkdiag(Z_2, \ldots, Z_{r+1})$, whose columns span $\ker \kronG$.
	Moreover, we introduce the matrix of perturbations $\kronE \coloneqq \blkdiag(E_2, \ldots, E_{r+1})$ such that $\kronZ = (\id_r \otimes \overline Z) + \varepsilon \, \kronE$.
	To show that $\kronZ^\transp \paren[big](){\kronD - \kronM \kronF(\xVec)} \kronZ$ is non-singular, we consider the first term:
	\begin{equation*}
		\begin{aligned}
			\kronZ^\transp \kronD \, \kronZ
			&
			=
			(\id_r \otimes \overline Z)^\transp \kronD\, (\id_r \otimes \overline Z)
			+
			\varepsilon \, \kronE^\transp \kronD (\id_r \otimes \overline Z) 
			+
			\varepsilon \, (\id_r \otimes \overline Z)^\transp \kronD \kronE 
			+
			\varepsilon^2 \, \kronE^\transp \kronD \kronE 
			\\
			&
			=
			(\id_r \otimes \overline Z^\transp) (\overline D \otimes \id_n) (\id_r \otimes \overline Z)
			+
			\cO(\varepsilon)
			\\
			&
			=
			\overline D \otimes (\overline Z^\transp \overline Z)
			+
			\cO(\varepsilon)
			.
		\end{aligned}
	\end{equation*}
	Here we utilized $(A \otimes B)(C \otimes D) = (AC) \otimes (BD)$. 
	By \Cref{lemma:Dbar_is_invertible}, $\overline D$ is nonsingular, thus $\overline D \otimes (\overline Z^\transp \overline Z)$ is nonsingular.
	The second term $\kronZ^\transp \kronM \kronF(\xVec) \kronZ$ is of order~$\cO(\Delta_0)$.
	Consequently, $\kronZ^\transp \paren[big](){\kronD - \kronM \kronF(\xVec)} \kronZ$ is nonsingular for sufficiently small $\varepsilon > 0$ and $\Delta_0 > 0$.
\end{proof}

The techniques used in the previous proof can also be applied to prove the solvability of system~\eqref{eq:cG_time_step_simplified}, which was claimed in~\Cref{thm:solvability}. 

\begin{proof}[Proof of~\Cref{thm:solvability}]
	We first note that the matrices~$D$ and~$M$ from~\eqref{eq:defDM} may also be defined as integrals over unit interval~$[0,1]$ with accordingly transformed Lagrange basis functions. This then shows that~$D$ is independent of~$\Delta$, whereas~$M$ (and thus $\kronM$) scales with~$\Delta$. More precisely, we have~$\kronM = \Delta \widehat \kronM$ with~$\widehat \kronM$ being independent of~$\Delta$. 
	Second, we rewrite~\eqref{eq:cG_time_step_simplified} in the form 	
	\begin{align*}
		\kronD \xVec - \Delta \widehat \kronM \fVec(\xVec) + \kronG(\xVec)^\transp \lVec 
		&
		= 
		\bVec
		, 
		\\
		\gVec(\xVec) 
		&
		= 
		0
	\end{align*}
	with $\fVec(\xVec) \coloneqq [f(x_1,t_1); \ldots; f(x_r,t_r)] \in \R^{r n}$, $\gVec(\xVec) \coloneqq [g(x_1,t_1); \ldots; g(x_r,t_r)] \in \R^{r m}$, and~$\bVec \coloneqq - D_1\otimes x_0 + M_1\otimes f(x_0,t_0) \in \R^{r n}$. 	
	We show the local solvability of this system by the implicit function theorem. 
	For this, we consider the function~$h\colon \R^{1+r(n+m)} \to \R^{r(n+m)}$, defined by
	\begin{equation*}
		h(\Delta,\xVec,\lVec)
		\coloneqq 
		\begin{bmatrix}
			\kronD \xVec - \Delta \widehat \kronM \fVec(\xVec) + \kronG(\xVec)^\transp \lVec - \bVec
			\\
			\gVec(\xVec)
		\end{bmatrix}
		.
	\end{equation*}
	Obviously, this function is sufficiently smooth by \cref{assumption:standing_assumption}. 
	Further, the Jacobian \wrt~$(\xVec,\lVec)$ equals~\eqref{eq:Newton_matrix} and it is invertible for sufficiently small~$\Delta$ by~\Cref{proposition:invertibility_Newton_matrix}. 
	To apply the implicit function theorem we finally need a root of~$h$ with~$\Delta=0$, \ie, we have to show the solvability of 
	\begin{subequations}
		\label{eq:inProof:Delta0System}
		\begin{align}
			\kronD \xVec + \kronG(\xVec)^\transp \lVec 
			&
			= 
			\bVec
			, 
			\\
			\gVec(\xVec) 
			&
			=
			0 
		\end{align}
	\end{subequations}
	for given~$\bVec$. 
	Due to the assumption on the structure of the constraint, \ie, $\gVec(\xVec) = \kronG \xVec + \tilde g(t)$, this follows again by~\Cref{proposition:invertibility_Newton_matrix} for sufficiently small $\Delta$.
\end{proof}

\begin{remark}
	In order to extend~\Cref{thm:solvability} to nonlinear constraints, one needs to ensure the solvability of system~\eqref{eq:inProof:Delta0System}. 
	This, however, may require additional assumptions on~$g$ and the initial data. 	
\end{remark}
%
%
%=============================================================================
%=========  Sect 3: Convergence
%=============================================================================
\section{Error Analysis}\label{sec:errAnalysis}
This section is devoted to the convergence analysis of the developed cG scheme~\eqref{eq:cG_time_step_simplified} for semi-explicit DAEs of the form~\eqref{eq:DAE}. Again we restrict the discussion to systems with constraints which are linear in~$x$. 
This means that we assume $g(x,t) = \tilde Bx - \tilde g(t)$ for a matrix~$\tilde B \in \R^{m\times n}$. 
We prove that the convergence rate for the cG scheme of degree $r$ equals~$r+1$. 
Since the approximation of the Lagrange multiplier is only defined as a linear combination of functionals, we consider as error measure the action of $\lambda-\Lambda$ on the constant function~$1$. 
The convergence proof is based on interpolation properties and corresponding convergence results for cG schemes applied to ODEs. 
%
%=============================================================================
\subsection{Convergence Results for ODEs}\label{sec:errAnalysis:ODE}
We summarize known convergence results for the cG approximation of an ODE, based on the variational formulation~\eqref{eq:variational_formulation_ODE}. 
Such schemes typically assume an exact integration so that we need to incorporate the quadrature error in a subsequent step.

% Step 1: no quadrature error 
We denote the idealized cG approximation on a subinterval $I=[0,\Delta]$ by 
\begin{equation*}
	\xt(t) = \sum_{j=0}^r \tilde x_j \, \varphi_j(t)
\end{equation*}
with coefficients $\tilde x_j \in \R^n$, $j = 0, \dots, r$. 
This means that, given $\tilde x_0$, $\xt$ satisfies 
\begin{equation}
	\label{eq:cGidealized}
	\int_I \big[ \dot {\xt}(t) - f(\xt(t),t) \big] \cdot \psi_i(t) \dt = 0
\end{equation}
for all $i=1,\dots,r$. 
Note that this does not include any quadrature error, since we have not replaced $f$ by its nodal interpolation. 
For this, the following convergence result holds. 

\begin{proposition}[\protect{see, \eg, \cite{EstepFrench:1994:1} or \cite[Sect.~5]{AkrivisMakridakis:2004:1}}]
	% EstF94: only r=1,2
	% AkrM04: parabolic case for arbitrary r
	\label{prop:ODEcase}
	Let $f$ be sufficiently smooth such that the solution of $\dot x = f(x,t)$ with initial value $x_0\in \R^n$ satisfies~$x\in C^{r+1}(0,T)$. 
	Further, let $\xt\in V$ be the idealized cG approximation defined by~\eqref{eq:cGidealized} of degree~$r$ and step size~$\Delta$ with $\xt(0)=x(0)=x_0$. 
	Then there exists a positive constant~$C$ such that 
	\begin{equation*}
		\norm{ x(t) - \xt(t)} 
		\le C\, \Delta^{r+1}\, \norm{ x }_{C^{r+1}(0,T)}
	\end{equation*}
	for all $t\in [0,T]$.
\end{proposition}
%
% Step 2: add quadrature error 
As outlined in~\Cref{sec:development:ODE}, the (computable) approximation in terms of the Lagrange basis calls for the application of a quadrature rule. 
The resulting approximation $\x(t) = \sum_{j=0}^r x_j \, \varphi_j(t)$ on an interval $I$ with coefficients $x_j \in \R^n$, $j=0, \dots, r$, is derived as the solution of 
\begin{equation}
	\label{eq:cGcomputed}
	\sum_{j=0}^r x_j \, D_{ij}
	-
	\sum_{j=0}^r f(x_j,t_j) \, M_{ij}
	=
	0
\end{equation}
for $i=1, \dots, r$.  
With $\kronD$ and $\kronM$ as introduced in~\Cref{sec:development}, the vector of unknown coefficients~$\xVec = [x_1; \ldots; x_r] \in \R^{rn}$, and $\fVec(\xVec) \coloneqq [f(x_1,t_1); \ldots; f(x_r,t_r)] \in \R^{r n}$, we can write~\eqref{eq:cGcomputed} as 
\begin{equation*}
	\kronD \xVec - \kronM \fVec(\xVec )
	= 
	\bVec 
	\coloneqq 
	- D_1\otimes x_0 + M_1\otimes f(x_0,t_0).
\end{equation*}
On the other hand, assuming identical initial data $x_0=\tilde x_0$, \eqref{eq:cGidealized} leads to the same system but with a perturbed right-hand side 
\begin{equation*}
	\tilde \bVec 
	\coloneqq 
	- D_1 \otimes x_0 + M_1 \otimes f(x_0,t_0) + \text{err}(\fVec). 
\end{equation*}
Here, $\text{err}(\fVec)$ denotes the quadrature error, which is of order $\Delta^{r+2}$ for a sufficiently smooth nonlinearity~$f$.  
Note, however, that the involved constant depends on the so-called Lebesgue constants 
\begin{align}
	\label{eq:LebesgueConstant}
	L_{r,\varphi}
	=
	L_{r,\varphi}(t_0, \dots, t_r)
	\coloneqq 
	\max_{t\in[0,\Delta]} \sum_{j=0}^r\, \abs{\varphi_j(t)}, \qquad
	L_{r,\psi}
	=
	L_{r,\psi}(t_1, \dots, t_r)
	\coloneqq 
	\max_{t\in[0,\Delta]} \sum_{i=1}^r \, \abs{\psi_i(t)}
	.
\end{align}
Both constants strongly depend on the choice of the Lagrange points. 
For equidistant~$t_j$ the Lebesgue constant grows exponentially in $r$, whereas Chebyshev nodes lead to a logarithmic growth only~\cite{Smith:2006:1}. 
Due to the invertibility of $\kronD$ and the fact that $\kronM$ scales with $\Delta$, the difference of the resulting coefficients $x_j$ and $\tilde x_j$ is again of order~$\Delta^{r+2}$ for sufficiently small step sizes. 

Considering the entire time interval $[0,T]$, these quadrature errors accumulate and lead to an additional factor of $\Delta^{-1}$. 
Note that~\Cref{assumption:left_and_right_Lagrange_points} implies that $x_0$ equals either the initial value or $x_r$ from the previous subinterval. 
Thus, errors in $x_r$ directly translate to perturbations of the right-hand side $b$. 
In total, this leads to the error estimate 
\begin{equation*}
	\norm{ x(t) - \x(t)} 
	\le C\, \Delta^{r+1} 
\end{equation*}
for all $t\in [0,T]$ under the assumptions of~\Cref{prop:ODEcase} with a constant~$C$ depending on $L_{r,\varphi}$ and $L_{r,\psi}$. 
%
%=============================================================================
\subsection{Convergence Results for Linear Constraints}
Based on the convergence results for ODEs, we analyze the cG scheme introduced in~\eqref{eq:cG_time_step_simplified}. 
Recall that we consider constraints that are linear in~$x$, \ie, we assume $g(x,t) = \tilde Bx - \tilde g(t)$. 
The full row rank property of $g_x$ from \Cref{assumption:standing_assumption} then translates into the full row rank (surjectivity) property of the matrix $\tilde B \in \R^{m\times n}$. 
Without loss of generality, we may assume that $\tilde B$ has the block structure~$\tilde B = [\, 0\hspace{5pt} B\,]$ with an invertible matrix~$B\in \R^{m\times m}$. This can be obtained by a simple transformation of variables based, \eg, on the QR decomposition of~$\tilde B$. 
Accordingly, we decompose the state variable into~$x = [y; z]$ with $y(t) \in \R^{n-m}$ and $z(t) \in \R^m$. Further, we denote by $f_1$ and $f_2$ the first $n-m$ and last $m$ components of $f$, respectively. These assumptions simplify system~\eqref{eq:DAE} to 
\begin{subequations}
	\label{eq:DAE:linear}
	\begin{align}
		\dot y & = f_1(y,z,t), \\
		\dot z & = f_2(y,z,t) - B^\transp\, \lambda, \\	
		0 & = B z - \tilde g(t).
	\end{align}
\end{subequations}

% cG
We apply the cG scheme \eqref{eq:cG_time_step_simplified} to the semi-explicit system~\eqref{eq:DAE:linear}. For this, we decompose the unknown coefficients $x_j$ in the same manner as the state variable, leading to  
\begin{equation*}
	\x(t) = \sum_{j=0}^r x_j \, \varphi_j(t), \qquad
	x_j = \begin{bmatrix} y_j \\ z_j\end{bmatrix}. 
\end{equation*}
Thus, we seek coefficients $y_j \in \R^{n-m}$ and $z_j \in \R^m$ for $j = 1, \dots, r$. 
Note that the resulting approximations of $y$ and $z$ (denoted by~$\y$ and $\z$, respectively) can also be written in terms of Lagrange polynomials. In this case, system~\eqref{eq:cG_time_step_simplified} turns into 
\begin{subequations}
	\label{eq:cG:linear}
	\begin{align}
		\sum_{j=0}^r \alert{y_j} \, D_{ij}
		- \sum_{j=0}^r f_1(\alert{y_j}, \alert{z_j}, t_j) \, M_{ij}
		\phantom{\ + B^\transp \alert{\lambda_i}}
		&= 0,  \label{eq:cG:linear_1} \\
		\sum_{j=0}^r \alert{z_j} \, D_{ij}
		- \sum_{j=0}^r f_2(\alert{y_j}, \alert{z_j}, t_j) \, M_{ij}
		+ B^\transp \alert{\lambda_i}
		&= 0,  \label{eq:cG:linear_2} \\
		B \alert{z_k}
		&
		= 
		\tilde g(t_k) 
		\label{eq:cG:linear_3}
	\end{align}
\end{subequations}
for $i, k = 1,\ldots, r$. 
Due to the invertibility of $B$, equation~\eqref{eq:cG:linear_3} determines all unknown coefficients $z_k$. 
With this, equation~\eqref{eq:cG:linear_1} can be solved for the coefficients~$y_j$ and, finally, equation~\eqref{eq:cG:linear_2} determines the coefficients~$\lambda_i$ of the Lagrange multiplier. 
Based on the convergence properties for ODEs, we derive the following result.
\begin{theorem}
	\label{thm:convergenceDAE}
	Consider a semi-explicit DAE~\eqref{eq:DAE} with sufficiently smooth right-hand side~$f$, a linear constraint $g(x,t) = \tilde Bx - \tilde g(t)$ with $\tilde B$ having full row rank, $\tilde g \in C^{r+2}(0,T)$, and a consistent initial condition, \ie, $\tilde B x(0) = \tilde g(0)$. 
	Further, we assume that~$x\in C^{r+1}(0,T)$ and Lagrange points with $t_0$ and $t_r$ being the starting and end point of the respective intervals. 
	Then, the cG approximation of degree~$r$ and step size~$\Delta$ given by~\eqref{eq:cG_time_step_simplified} satisfies
	\begin{align*}
		\norm{x(t)-\x(t)} 
		&
		\le
		C\, \Delta^{r+1}
		, 
		\\
		\norm{\lambda - \Lambda}_{*,I}
		\coloneqq
		\sup_{0 \neq v \in C^\infty(I)} \frac{ \dual{\lambda - \Lambda}{v} }{ \norm{v}_{C^\infty(I)} }
		&
		\le
		C\, \Delta^{r+1}
	\end{align*}
	for all $t \in [0,T]$ and all sub-intervals~$I = [(\ell-1)\Delta, \ell\Delta]$, $\ell = 1, \dots, T/\Delta$. 
	Therein, the constant~$C$ includes the Lebesgue constants~$L_{r,\varphi}$ and $L_{r,\psi}$,  which in turn depend on the distribution of the Lagrange points~$t_j$. 
\end{theorem}
\begin{proof}
	Without loss of generality, we assume that the DAE is of the form~\eqref{eq:DAE:linear} and we consider a single time step on the interval $I = [0, \Delta]$. 
	Since the initial data is consistent and~$t_r$ equals the end point of the interval~$I$, equation~\eqref{eq:cG:linear_3} guarantees that we obtain consistent initial values in every time step. 
	More precisely, this means that the prescribed coefficient $x_0 = [y_0; z_0]$ at time $t_0 = 0$ satisfies $B z_0 = \tilde B x_0 = \tilde g(0)$. 

	% z
	Considering equation~\eqref{eq:cG:linear_3} for $k=1,\dots,r$, we obtain 
	\begin{equation*}
		\z(t) 
		= 
		\sum_{j=0}^r z_j \, \varphi_j(t)
		= 
		B^{-1} \sum_{j=0}^r \tilde g(t_j) \, \varphi_j(t)
		,
	\end{equation*}
	\ie, the resulting approximation~$\z$ equals the interpolation polynomial of $B^{-1}\tilde g$ in the Lagrange points $t_0,\dots, t_r$. 
	This directly implies the error estimate
	\begin{align}
		\label{eq:proofZ}
		\norm{z(t) - \z(t)}
		\le 
		C_B \, \norm{\tilde g(t) - \Pi_r\tilde g(t)}
		\le 
		\tfrac{1}{(r+1)!}\, C_B\, \Delta^{r+1} \norm{\tilde g}_{C^{r+1}(0,T)}
		.
	\end{align}
	Note that this result may be improved by an appropriate choice of the interpolation points. Further, it holds $z(t_j) = \z(t_j)$ in all Lagrange points. 

	% y
	Next, we consider the remaining part of the state vector, namely the $y$ component. For this, we consider equation~\eqref{eq:cG:linear_1} and note that this is nothing else than the standard cG scheme of degree~$r$ applied to the ODE
	\begin{equation*} 
		\dot y 
		= \tilde f_1(y,t) 
		\coloneqq f_1(y, B^{-1}\tilde g,t).
	\end{equation*}
	Thus, with the assumed regularity of~$f$ and~$\tilde g$, the error analysis for the ODE case in~\Cref{sec:errAnalysis:ODE} implies that the error $y(t)-\y(t)$ is of order~$r+1$ as well. Together with~\eqref{eq:proofZ} this yields the stated error bound for the state variable~$x$.

	% Lagrange
	Finally, we consider the Lagrange multiplier. The approximation $\Lambda = \sum_{k=1}^r \lambda_k \, \xi_k$ is defined by equation~\eqref{eq:cG:linear_2}, whereas the exact multiplier satisfies 
	\begin{equation*}
		B^\transp\, \lambda(t)
		= 
		f_2(y,z,t) - \dot z(t)
		= 
		f_2(y,B^{-1}\tilde g,t) - B^{-1}\dot {\tilde g}(t)
	\end{equation*}
	due to $Bz(t) = \tilde g(t)$.
	For an estimate of the action of the difference $\lambda - \Lambda$ on an arbitrary smooth test function~$v\in C^\infty(I)$, we consider 
	% [use embedding such that L2 inner product appears for smooth functions] 
	%
	\begin{align*}
		\norm[big]{B^\transp \dual{\lambda - \Lambda}{v}}
		&
		= 
		\norm[Big]{\int_0^\Delta B^\transp\lambda(t)v(t) \dt - \sum_{k=1}^r B^\transp\lambda_k v(s_k)} 
		\\
		& 
		\le 
		\norm[Big]{\int_0^\Delta B^{-1}\dot {\tilde g}(t)\, v(t) \dt - \sum_{k=1}^r\sum_{j=0}^r z_j D_{kj} v(s_k)} 
		\\
		&
		\qquad
		+ \norm[Big]{\int_0^\Delta f_2(y,B^{-1}\tilde g,t)\, v(t) \dt - \sum_{k=1}^r\sum_{j=0}^r f_2(y_j, z_j, t_j) \, M_{kj} v(s_k)} 
		. 
	\end{align*}
	For the first term on the right-hand side, we use the Lagrange interpolation~$v(t) \approx \sum_{k=1}^r \psi_k(t)\, v(s_k)$ and estimate~\eqref{eq:proofZ}. This then leads to 
	\begin{align*}
		\sum_{k=1}^r\sum_{j=0}^r z_j D_{kj} v(s_k)
		= 
		\int_0^\Delta \Big( \sum_{j=0}^r z_j\, \dot \varphi_j(t)\Big) \, \Big(\sum_{k=1}^r \psi_k(t)\, v(s_k)\Big) \dt
		&
		=
		\int_0^\Delta \dot \z(t) \, v(t) \dt + E_\text{int}
		\\
		&
		=
		\int_0^\Delta B^{-1}\dot {\tilde g}(t)\, v(t) \dt + E_\text{int} + E_z 
	\end{align*}
	with error terms being bounded by  
	\begin{equation*}
		\norm{E_\text{int}}
		%	=  
		%	\int_0^\Delta \dot \z(t)\, \Big(\sum_{k=1}^r \psi_k(t)\, v(s_k) - v(t)\Big) \dt	
		%	\le 
		%	C\, \tfrac{1}{r!} \, \Delta^r \norm{v}_{C^r(0,\Delta)} 
		%	\int_0^\Delta \norm{ \dot \z(t) } \dt 
		\le 
		C\, \Delta^{r+1}\, \norm{v}_{C^r(I)}, \qquad  
		\norm{E_z}	
		%	=  
		%	\int_0^\Delta \dot \z(t) \, v(t) \dt - \int_0^\Delta \dot z(t) \, v(t) \dt 
		\le 
		C\, \Delta^{r+1}\, \norm{v}_{C^0(I)}.
	\end{equation*}
	For the second part, we proceed similarly, leading to 
	\begin{align*}
		\sum_{k=1}^r\sum_{j=0}^r f_2(y_j, z_j, t_j) \, M_{kj} v(s_k)
		&
		= 
		\int_0^\Delta \Big( \sum_{j=0}^r f_2(y_j, B^{-1} {\tilde g}(t_j), t_j)\, \varphi_j(t) \Big)\, \Big( \sum_{k=1}^r \psi_k(t)\, v(s_k) \Big) \dt
		\\
		&
		= 
		\int_0^\Delta f_2(y(t), B^{-1} {\tilde g}(t), t)\, \Big( \sum_{k=1}^r \psi_k(t)\, v(s_k) \Big) \dt + E_f  
		\\
		&
		= 
		\int_0^\Delta f_2(y(t), B^{-1} {\tilde g}(t), t)\, v(t) \dt
		+ E_f + E_\text{int}. 
	\end{align*}
	Using the local Lipschitz continuity of $f$, the error estimate of~$y-\y$, and the Lebesgue constants~$L_{r,\varphi}$, $L_{r,\psi}$ defined in~\eqref{eq:LebesgueConstant}, the error terms can be bounded by 
	\begin{equation*}
		\norm{E_f}
		%	= 
		%	\int_0^\Delta \Big( \sum_{j=0}^r f_2(y_j, B^{-1} {\tilde g}(t_j), t_j)\, \varphi_j(t) - f_2(y(t), B^{-1} {\tilde g}(t), t) \Big)\, \Big( \sum_{k=1}^r \psi_k(t)\, v(s_k) \Big) \dt
		\le 
		C\, \Delta^{r+2}\, L_{r,\varphi}\, L_{r,\psi}\, \norm{v}_{C^0(I)}, \qquad
		\norm{E_\text{int}}
		\le 
		C\, \Delta^{r+1}\, \norm{v}_{C^r(I)}
		.
	\end{equation*}
	Altogether, we obtain the claimed error estimate of order~$\Delta^{r+1}$ in the dual norm. 
\end{proof}
\begin{remark}
	\label{rem:LagrangeConstTestfct}
	Considering the constant function~$1_I$ (indicator function) as test function, we can show with the help of \Cref{lemma:row_column_sums_of_Dbar} that~$\norm{\dual{\lambda - \Lambda}{1_I}} \le C\, \Delta^{r+2}$ for every sub-interval~$I = [(\ell-1) \Delta, \ell\Delta]$, $\ell = 1, \dots, T/\Delta$. 
	Hence, for this particular test function, the error is of one order higher than the error for the state. 
\end{remark}
To summarize, the proposed cG scheme maintains the convergence rate from the ODE case and this rate even extends to the Lagrange multiplier. 
We close the paper with a numerical illustration of the convergence results. 
%
%
%=============================================================================
%=========  Sect 4: Numerics
%=============================================================================
\section{Numerical Experiments}
\label{sec:numerics}
This final section is devoted to three numerical experiments, which verify the theoretical findings but also show the applicability to related DAE models of higher index (not currently covered by our theory). 
We would like to emphasize once more that the cG scheme~\eqref{eq:cG_time_step_simplified} is directly applied to the DAE without any index reduction or regularization. 
Nevertheless, no artificial drift-off can be observed since the constraint is explicitly enforced in the Lagrange points. 
Thus, the error in the constraint is directly linked to the accuracy of the nonlinear solver applied to solve \eqref{eq:cG_time_step_simplified}. 
In our implementation, we currently use \matlab's \verb!fsolve! for this purpose, with \verb!OptimalityTolerance! set to $10^{-13}$.

In all experiments, we use cG schemes of various degrees up to $r = 5$ and with equidistant Lagrange points $t_0, \ldots, t_r$. We contend ourselves with subintervals of constant length~$\Delta$, leaving the discussion of adaptivity to future work. 

The code for this paper is available on GitHub\footnote{https://github.com/rolandherzog/cg-schemes-for-daes} and can be cited as \cite{AltmannHerzog:2021:1}.
%
%=============================================================================
\subsection{Simple Circuit Example}
\label{subsection:simple_linear_circuit}
% https://depositonce.tu-berlin.de/bitstream/11303/1858/1/Dokument_18.pdf

We follow an example presented in~\cite[Sect.~5.2.1]{Baechle:2007:1}.
Eliminating the trivial variables in this example, we obtain the linear system
\begin{align*}
	\dot q_1 & = - \sin(100\,t) - i_V, \\
	\dot q_2 & = - q_2 - \sin(100\,t) - i_V, \\  
	0 & = q_1 + q_2 - \sin(100\,t)
\end{align*}
with unknowns~$q_1$, $q_2$, and~$i_V$. 
In our setting, $i_V$ plays the role of a Lagrange multiplier. 
Starting with initial data~$q_1(0) = q_2(0) = 0$, we can conclude $i_V(0)=-50$ by consistency of the system. 
For this system, the exact solution is explicitly known so we can use it to evaluate the discretization error.
As time horizon we set~$T = 1$. 

% convergence 
Recall that we have shown in \cref{thm:convergenceDAE} that the cG scheme of degree~$r$ converges with order~$r+1$. 
In the convergence history shown in \cref{fig:circuit:states}, we observe the rates~$r+1$, and~$r+2$ for even polynomial degrees. 
Note that this increase of the rate results from the choice of equidistant Lagrange points. For Gau\ss-Lobatto points, one can observe the rates~$2r$.
As a comparison, we also plot the convergence history of the well-known Radau~IIa methods of order~$3$ and~$5$. 
As an error measure for the Lagrange multiplier, we have considered the application of the difference~$\lambda - \Lambda$ to the constant function on a subinterval, \ie, $\norm{\dual{\lambda - \Lambda}{1_I}}$. Here, we always consider the last subinterval~$I=[T-\Delta, T]$. 
The results are again in agreement with \cref{rem:LagrangeConstTestfct}, see \cref{fig:circuit:lambda}. 
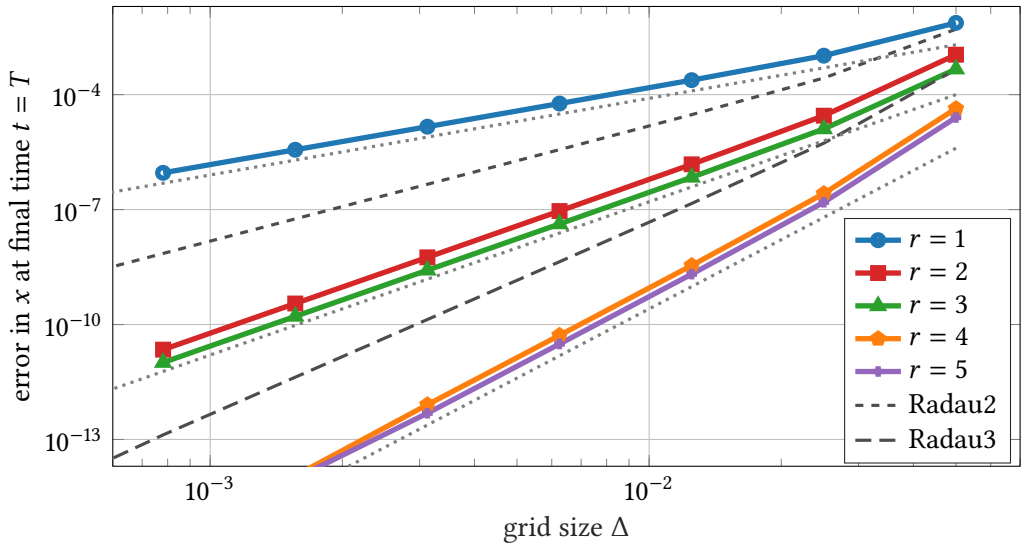
\begin{figure}
	\centering

% color for tikz images
\definecolor{color0}{rgb}{0.12156862745098,0.466666666666667,0.705882352941177}
\definecolor{color1}{rgb}{1,0.498039215686275,0.0549019607843137}
\definecolor{color2}{rgb}{0.172549019607843,0.627450980392157,0.172549019607843}
\definecolor{color3}{rgb}{0.83921568627451,0.152941176470588,0.156862745098039}
\definecolor{color4}{rgb}{0.580392156862745,0.403921568627451,0.741176470588235}

\begin{tikzpicture}

\begin{axis}[%
width=4.7in,
height=2.4in,
scale only axis,
xmode=log,
xmin=0.0006,
xmax=0.07,
xminorticks=true,
xlabel style={font=\color{white!15!black}},
xlabel={grid size $\Delta$},
ylabel={error in $x$ at final time~$t=T$},
ymode=log,
ymin=2e-14,
ymax=0.02,
yminorticks=true,
title style={font=\bfseries},
xmajorgrids,
%xminorgrids,
ymajorgrids,
%yminorgrids,
legend style={legend cell align=left, align=left, draw=white!15!black, at={(0.99,0.54)}}
]
\addplot [color=color0, line width=2.0pt, mark=o]
table[row sep=crcr]{%
%0.1	0.00433451979781779\\
0.05	0.00746764151256704\\
0.025	0.00104014780757267\\
0.0125	0.000238883402168628\\
0.00625	5.85830808079364e-05\\
0.003125	1.45771034981934e-05\\
0.0015625	3.64002095932221e-06\\
0.00078125	9.09739870732153e-07\\
};
\addlegendentry{$r=1$}

\addplot [color=color3, line width=2.0pt, mark=square*]
  table[row sep=crcr]{%
%0.1	0.00851184251212183\\
0.05	0.00110226298571773\\
0.025	2.82018683414164e-05\\
0.0125	1.516864199175e-06\\
0.00625	9.15457001652455e-08\\
0.003125	5.67261593367985e-09\\
0.0015625	3.53793683244905e-10\\
0.00078125	2.21010139715475e-11\\
};
\addlegendentry{$r=2$}

\addplot [color=color2, line width=2.0pt, mark=triangle*]
table[row sep=crcr]{%
%0.1	0.00187950082021352\\
0.05	0.000463995925865373\\
0.025	1.27412683597476e-05\\
0.0125	6.96014942753954e-07\\
0.00625	4.21637934589882e-08\\
0.003125	2.6151069823861e-09\\
0.0015625	1.6313543262367e-10\\
0.00078125	1.01884869736449e-11\\
};
\addlegendentry{$r=3$}

\addplot [color=color1, line width=2.0pt, mark=pentagon*]
table[row sep=crcr]{%
0.05	4.34110948173315e-05\\
0.025	2.62191444285734e-07\\
0.0125	3.47627333524517e-09\\
0.00625	5.22674751464339e-11\\
0.003125	8.07655560746934e-13\\
0.0015625	1.1226161079589e-14\\
0.00078125	4.47476350724876e-15\\
};
\addlegendentry{$r=4$}

\addplot [color=color4, line width=2.0pt, mark=+]
table[row sep=crcr]{%
	0.05	2.53269893358441e-05\\
	0.025	1.55346469595784e-07\\
	0.0125	2.06791362783479e-09\\
	0.00625	3.11216269221367e-11\\
	0.003125	4.84609037372747e-13\\
	0.0015625	9.49905937503684e-15\\
	0.00078125	1.0833637964918e-14\\
};
\addlegendentry{$r=5$}

%%%%%%%% alternative integration schemes %%%%%%%%

%\addplot [color=red]
%table[row sep=crcr]{%
%0.1	0.0115222246685953\\
%0.05	0.000626750518172585\\
%0.025	0.00296120852331292\\
%0.0125	0.00176345618371371\\
%0.00625	0.000943629626992233\\
%0.003125	0.000486856596248898\\
%0.0015625	0.000247166696252947\\
%0.00078125	0.000124517081463872\\
%0.000390625	6.24919822614517e-05\\
%};
%\addlegendentry{impl Euler}

%\addplot [color=blue, dashed]
%table[row sep=crcr]{%
%0.1	6.04845298104123\\
%0.05	2.21181688856449\\
%0.025	0.0290026798760841\\
%0.0125	0.113179842821305\\
%0.00625	0.0222369013322102\\
%0.003125	0.00496332857765357\\
%0.0015625	0.00118837947914266\\
%0.00078125	0.000292168791419823\\
%0.000390625	7.25362604941877e-05\\
%%0.0001953125	1.81025345719727e-05\\
%%9.765625e-05	4.52404004799219e-06\\
%%4.8828125e-05	1.13266187657134e-06\\
%%2.44140625e-05	2.79731848488604e-07\\
%};
%\addlegendentry{ode23}

%\addplot [color=black!40!gray, line width=1.2pt, dashed]
%table[row sep=crcr]{%
%0.1	0.0115222246685954\\
%0.05	0.000626750518172507\\
%0.025	0.00296120852331296\\
%0.0125	0.0017634561837136\\
%0.00625	0.000943629626992115\\
%0.003125	0.000486856596248898\\
%0.0015625	0.000247166696252947\\
%0.00078125	0.000124517081463675\\
%0.000390625	6.24919822608237e-05\\
%0.0001953125	3.1304359795794e-05\\
%9.765625e-05	1.56667735520081e-05\\
%};
%\addlegendentry{Radau IIa $s=1$}

\addplot [color=black!40!gray, line width=1.2pt, dashed]
table[row sep=crcr]{%
%0.1	0.0148721631351414\\
0.05	0.00510276063376635\\
0.025	0.000271426200869743\\
0.0125	3.00121969695759e-05\\
0.00625	3.66833535128702e-06\\
0.003125	4.57302841420496e-07\\
0.0015625	5.72031598631257e-08\\
0.00078125	7.15656456307633e-09\\
0.000390625	8.95071321935005e-10\\
0.0001953125	1.11917878304379e-10\\
9.765625e-05	1.3994744364793e-11\\
};
%\addlegendentry{Radau IIa $s=2$}
\addlegendentry{Radau$2$}

\addplot [color=black!40!gray, line width=1.2pt, dashed, dash pattern=on 6pt off 3pt]
table[row sep=crcr]{%
%	0.1	0.00814877934289687\\
	0.05	0.000455449006945745\\
	0.025	5.40504419193982e-06\\
	0.0125	1.44285330140141e-07\\
	0.00625	4.35666233007591e-09\\
	0.003125	1.3515744482234e-10\\
	0.0015625	4.21832815661686e-12\\
	0.00078125	1.3137749234417e-13\\
	0.000390625	3.29719416323593e-15\\
	0.0001953125	2.2377260456559e-15\\
	9.765625e-05	5.49532360539321e-15\\
};
%\addlegendentry{Radau IIa $s=3$}
\addlegendentry{Radau$3$}

%%%%%%%% order lines %%%%%%%%

\addplot [color=gray, line width=1.2pt, dotted]
table[row sep=crcr]{%
	0.05	0.002\\
	0.0005	0.0000002\\
};
%\addlegendentry{order $2$}

\addplot [color=gray, line width=1.2pt, dotted]
table[row sep=crcr]{%
	0.05	0.0001\\
	0.0005	0.000000000001\\
};
%\addlegendentry{order $4$}

\addplot [color=gray, line width=1.2pt, dotted]
table[row sep=crcr]{%
	0.05		0.000004\\
	0.0005	    0.000000000000000004\\
};
%\addlegendentry{order $6$}

\end{axis}
\end{tikzpicture}%
	\caption{Convergence orders for the state obtained by the cG scheme applied to the circuit example from \cref{subsection:simple_linear_circuit}. For comparison, also the results of Radau~IIa with two and three stages are shown, which converge with order~$3$ and~$5$, respectively. The gray dotted lines indicate the orders $2$, $4$, and~$6$.}
	\label{fig:circuit:states}
\end{figure}
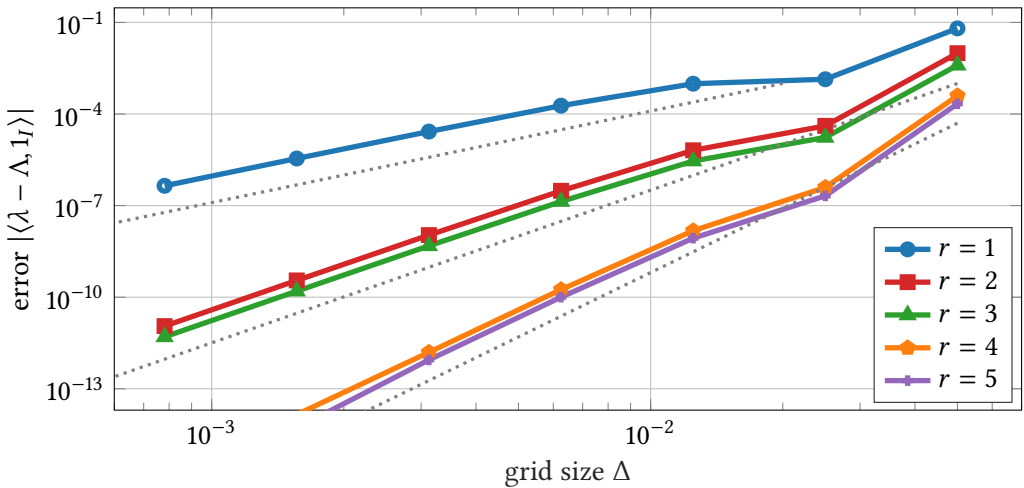
\begin{figure}
	\centering

% color for tikz images
\definecolor{color0}{rgb}{0.12156862745098,0.466666666666667,0.705882352941177}
\definecolor{color1}{rgb}{1,0.498039215686275,0.0549019607843137}
\definecolor{color2}{rgb}{0.172549019607843,0.627450980392157,0.172549019607843}
\definecolor{color3}{rgb}{0.83921568627451,0.152941176470588,0.156862745098039}
\definecolor{color4}{rgb}{0.580392156862745,0.403921568627451,0.741176470588235}

\begin{tikzpicture}

\begin{axis}[%
width=4.7in,
height=2.1in,
scale only axis,
xmode=log,
xmin=0.0006,
xmax=0.07,
xminorticks=true,
xlabel style={font=\color{white!15!black}},
xlabel={grid size $\Delta$},
ylabel={error $|\langle \lambda-\Lambda,1_I\rangle|$},
ymode=log,
ymin=2e-14,
ymax=0.3,
yminorticks=true,
title style={font=\bfseries},
xmajorgrids,
%xminorgrids,
ymajorgrids,
%yminorgrids,
legend style={legend cell align=left, align=left, draw=white!15!black, at={(0.99,0.455)}}
]
\addplot [color=color0, line width=2.0pt, mark=o]
table[row sep=crcr]{%
%0.1	0.0389560187549958\\
0.05	0.0637276930892237\\
0.025	0.00137572921719253\\
0.0125	0.000982439701667337\\
0.00625	0.000187280391966055\\
0.003125	2.65896735690135e-05\\
0.0015625	3.48810987151754e-06\\
0.00078125	4.45119239470687e-07\\
};
\addlegendentry{$r=1$}

\addplot [color=color3, line width=2.0pt, mark=square*]
table[row sep=crcr]{%
%0.1	0.0782782432013865\\
0.05	0.00982233403180355\\
0.025	4.09100647076333e-05\\
0.0125	6.47409217613415e-06\\
0.00625	3.05024779989616e-07\\
0.003125	1.07985551750511e-08\\
0.0015625	3.53994153434734e-10\\
0.00078125	1.12934279983268e-11\\
};
\addlegendentry{$r=2$}

\addplot [color=color2, line width=2.0pt, mark=triangle*]
table[row sep=crcr]{%
%0.1	0.0154870403870535\\
0.05	0.00395853226016567\\
0.025	1.71010770699809e-05\\
0.0125	2.88228608102958e-06\\
0.00625	1.35875761753113e-07\\
0.003125	4.81013719133561e-09\\
0.0015625	1.57654868326862e-10\\
0.00078125	5.02860947326766e-12\\
};
\addlegendentry{$r=3$}

\addplot [color=color1, line width=2.0pt, mark=pentagon*]
table[row sep=crcr]{%
%0.1000    0.016510531796884\\
0.05	0.000395482973292588\\
0.025	3.95428369381889e-07\\
0.0125	1.50723175051937e-08\\
0.00625	1.77196118889e-10\\
0.003125	1.56837043352454e-12\\
0.0015625	1.43357548054723e-14\\
0.00078125	1.1518563880486e-15\\
};
\addlegendentry{$r=4$}

\addplot [color=color4, line width=2.0pt, mark=+]
table[row sep=crcr]{%
%0.1000    0.008082820966216\\
0.05	0.000216626102645875\\
0.025	2.08693610836974e-07\\
0.0125	8.56418957972949e-09\\
0.00625	1.00295494132041e-10\\
0.003125	8.84237127962706e-13\\
0.0015625	5.22498710964214e-15\\
0.00078125	2.17187379192296e-15\\
};
\addlegendentry{$r=5$}

\addplot [color=gray, line width=1.2pt, dotted]
table[row sep=crcr]{%
	0.02	0.001\\
	0.0002	0.000000001\\
};
%\addlegendentry{order $3$}

\addplot [color=gray, line width=1.2pt, dotted]
table[row sep=crcr]{%
	0.05	0.001\\
	0.0005	0.0000000000001\\
};
%\addlegendentry{order $5$}

\addplot [color=gray, line width=1.2pt, dotted]
table[row sep=crcr]{%
	0.05	0.00005\\
	0.0005	0.0000000000000000005\\
};
%\addlegendentry{order $7$}

\end{axis}
\end{tikzpicture}%
 % time interval I
	\caption{Convergence orders for the Lagrange multipliers on the final subinterval obtained by the cG scheme applied to the circuit example from \cref{subsection:simple_linear_circuit}. Here, the gray dotted lines indicate orders $3$, $5$, and $7$.}
	\label{fig:circuit:lambda}	
\end{figure}
%
%
%=============================================================================
\subsection{Quasilinear Heat Equation}
\label{subsection:quasilinear_heat_equation}
In this second example, we consider a coupling of two one-dimensional quasilinear heat equations. 
Both problems are coupled through a transmission condition representing a thermal resistance.
To be more concrete, we consider the problem
\begin{subequations}
\label{eq:quasilinear_heat_equations}
\begin{align}
		\dot u - [u(z)^{c_1}]'' 
		&
		= 
		0
		\qquad \text{in } (0,1)
		, 
		\\
		\dot u - [u(z)^{c_2}]'' 
		&
		= 
		0 
		\qquad \text{in } (1,2)
		.
\end{align}
\end{subequations}
Here $[\, \cdot\, ]'$ denotes differentiation with respect to the spatial variable~$z$.
Problem \eqref{eq:quasilinear_heat_equations} is discretized using a standard finite difference scheme described in \cite{Gipson:1987:1}.
Specifically, suppose that \verb!K! denotes the standard stiffness matrix with natural (Neumann) boundary conditions on either interval,
\begin{equation*}
	\verb!K!
	=
	\frac{1}{h^2}
	\begin{bmatrix}
		1 & -1
		\\
		-1 & 2 & -1 
		\\
		& \ddots & \ddots & \ddots 
		\\
		& & -1 & 2 & -1 
		\\
		& & & -1 & 1
	\end{bmatrix}
	.
\end{equation*}
In our experiments we use an equidistant spatial grid size of $h = 1/40$ on the spatial domain $[0,1] \cup [1,2]$, \ie, the discrete state is of dimension~$n = 82$.
We denote by \lstinline!x! the vector approximating the values of this state.
Then the right hand side~$f$ in \matlab notation is given by
\begin{lstlisting}[frame=none,backgroundcolor=]
	f = - [K*x(1:end/2).^c1; K*x(end/2+1:end).^c2];
\end{lstlisting}

\begin{figure}
	\centering
	\input{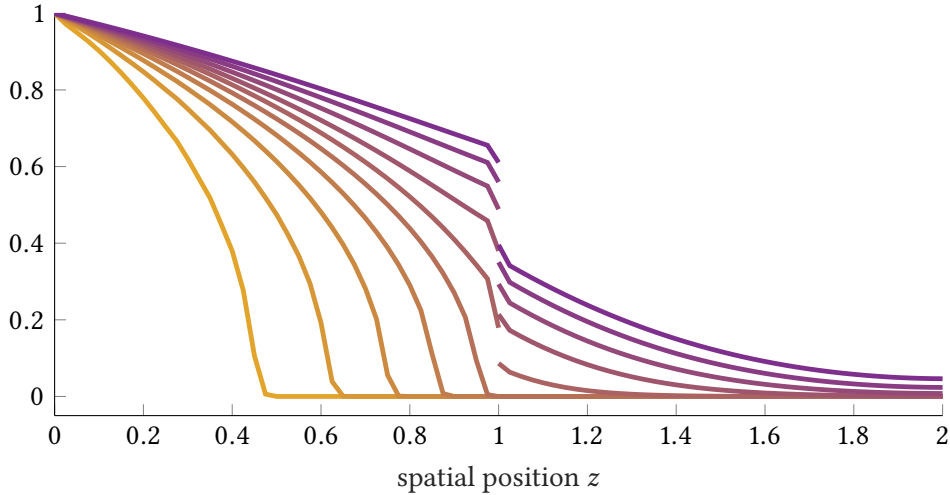}
	\caption{Evolution of the temperature obtained by applying the cG scheme with $r = 1$ to the quasilinear heat equation example from \cref{subsection:quasilinear_heat_equation} with $c_1 = 3$ and $c_2 = 1$. Each line indicates one point in time~$t = 4 \ell \Delta$, $\ell = 1, \dots, 10$ (from orange to purple).} 
	\label{fig:resistence:state}
\end{figure}

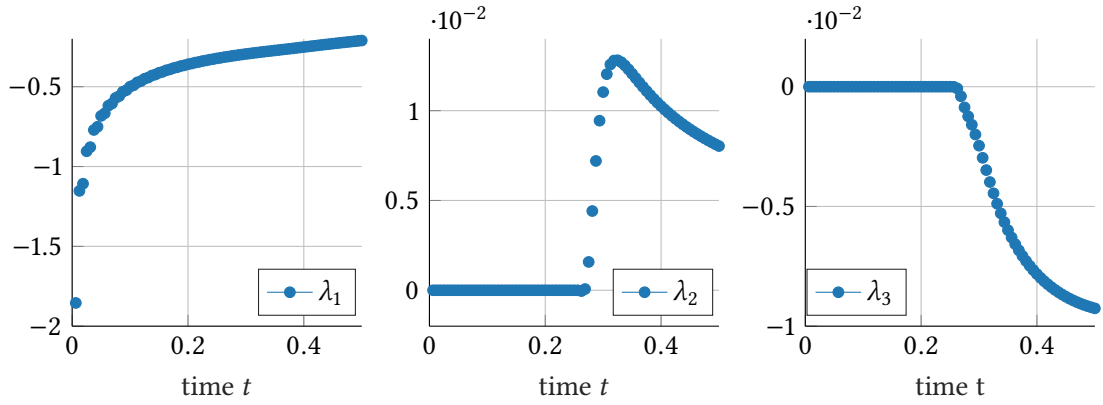
\begin{figure}
	\centering

% color for tikz images
\definecolor{color0}{rgb}{0.12156862745098,0.466666666666667,0.705882352941177}
\definecolor{color1}{rgb}{1,0.498039215686275,0.0549019607843137}
\definecolor{color2}{rgb}{0.172549019607843,0.627450980392157,0.172549019607843}
\definecolor{color3}{rgb}{0.83921568627451,0.152941176470588,0.156862745098039}
\definecolor{color4}{rgb}{0.580392156862745,0.403921568627451,0.741176470588235}

\begin{tikzpicture}

\begin{axis}[%
width=1.5in,
height=1.5in,
at={(0.758in,0.481in)},
scale only axis,
xmin=0,
xmax=0.5,
xlabel style={font=\color{white!15!black}},
xlabel={time $t$},
ymin=-2,
ymax=-0.2,
axis background/.style={fill=white},
title style={font=\bfseries},
axis x line*=bottom,
axis y line*=left,
xmajorgrids,
ymajorgrids,
legend style={legend cell align=left, align=left, draw=white!15!black, at={(0.98,0.20)}}
]
\addplot [color=color0, draw=none, mark=*]
  table[row sep=crcr]{%
0.00625	-1.85455184020581\\
0.0125	-1.15258645022816\\
0.01875	-1.10751866506673\\
0.025	-0.904587843445093\\
0.03125	-0.878534162977043\\
0.0375	-0.770806000334235\\
0.04375	-0.750843573947501\\
0.05	-0.682735009817449\\
0.05625	-0.666544007408652\\
0.0625	-0.619117665890341\\
0.06875	-0.605559142932262\\
0.075	-0.570409668115229\\
0.08125	-0.55880838043585\\
0.0875	-0.531581797729053\\
0.09375	-0.521493483221748\\
0.1	-0.499696680963364\\
0.10625	-0.490809412363616\\
0.1125	-0.47290792847636\\
0.11875	-0.464995988831664\\
0.125	-0.449992778198947\\
0.13125	-0.442885910737727\\
0.1375	-0.430101368579627\\
0.14375	-0.423668304190625\\
0.15	-0.412624015301883\\
0.15625	-0.406762115867802\\
0.1625	-0.397109944163091\\
0.16875	-0.39173729059122\\
0.175	-0.383217926042862\\
0.18125	-0.378268376014687\\
0.1875	-0.370684180837009\\
0.19375	-0.366103699804856\\
0.2	-0.359301284727551\\
0.20625	-0.355045189726924\\
0.2125	-0.348903526200272\\
0.21875	-0.344934532155391\\
0.225	-0.339356744799474\\
0.23125	-0.335643451651864\\
0.2375	-0.330551015865189\\
0.24375	-0.327066753967554\\
0.25	-0.322395298364242\\
0.25625	-0.319117241590273\\
0.2625	-0.314813463123841\\
0.26875	-0.31172193564915\\
0.275	-0.307741295974768\\
0.28125	-0.304819220770728\\
0.2875	-0.30112415957986\\
0.29375	-0.298356398291451\\
0.3	-0.294914148292615\\
0.30625	-0.292284902379274\\
0.3125	-0.289062080767583\\
0.31875	-0.286548399830281\\
0.325	-0.283503758268268\\
0.33125	-0.281070589462854\\
0.3375	-0.278152356601382\\
0.34375	-0.275753585607236\\
0.35	-0.272904619187987\\
0.35625	-0.270494080441508\\
0.3625	-0.267666788426125\\
0.36875	-0.265213759611098\\
0.375	-0.262381570199821\\
0.38125	-0.259876561642316\\
0.3875	-0.257033795760238\\
0.39375	-0.254484387623915\\
0.4	-0.251639585740104\\
0.40625	-0.249063111149506\\
0.4125	-0.246231871786337\\
0.41875	-0.243649247694715\\
0.425	-0.240849070920325\\
0.43125	-0.238280871615533\\
0.4375	-0.235528175906439\\
0.44375	-0.232992601973351\\
0.45	-0.230301325909439\\
0.45625	-0.227813419804175\\
0.4625	-0.225194589287917\\
0.46875	-0.222766157669945\\
0.475	-0.220227968044093\\
0.48125	-0.217867851689329\\
0.4875	-0.215415983754755\\
0.49375	-0.213130463128753\\
0.5	-0.210768474798879\\
};
\addlegendentry{$\lambda_1$}

\end{axis}
\end{tikzpicture}
%
%%% lam 2 %%%
%
\begin{tikzpicture}
\begin{axis}[%
width=1.5in,
height=1.5in,
at={(0.758in,0.481in)},
scale only axis,
xmin=0,
xmax=0.5,
xlabel style={font=\color{white!15!black}},
xlabel={time $t$},
ymin=-0.002,
ymax=0.014,
axis background/.style={fill=white},
title style={font=\bfseries},
axis x line*=bottom,
axis y line*=left,
xmajorgrids,
ymajorgrids,
legend style={legend cell align=left, align=left, draw=white!15!black, at={(0.98,0.20)}}
]
\addplot [color=color0, draw=none, mark=*]
table[row sep=crcr]{%
	0.00625	0\\
	0.0125	0\\
	0.01875	0\\
	0.025	0\\
	0.03125	0\\
	0.0375	0\\
	0.04375	2.2427976595239e-317\\
	0.05	1.43865311302196e-302\\
	0.05625	6.78454551522822e-289\\
	0.0625	2.04722680270872e-276\\
	0.06875	-1.47362199286149e-263\\
	0.075	1.74738968729016e-251\\
	0.08125	3.90969736457441e-238\\
	0.0875	9.84199987189839e-227\\
	0.09375	2.23374417276952e-217\\
	0.1	-2.36478528942403e-205\\
	0.10625	-8.96769269221498e-196\\
	0.1125	4.84723335848277e-185\\
	0.11875	-7.64688862733841e-175\\
	0.125	1.11740321521048e-166\\
	0.13125	1.78460756230906e-156\\
	0.1375	-2.9910140739454e-146\\
	0.14375	-3.09438266061278e-138\\
	0.15	-3.3914826195714e-129\\
	0.15625	4.14376903709e-120\\
	0.1625	-2.40625740806152e-113\\
	0.16875	4.76757760890188e-103\\
	0.175	-8.70514908193496e-96\\
	0.18125	5.75194552752922e-87\\
	0.1875	-1.47025294990972e-79\\
	0.19375	1.81798249535923e-71\\
	0.2	-4.89867346082086e-64\\
	0.20625	2.13044211356182e-56\\
	0.2125	-9.31587037030756e-49\\
	0.21875	9.31424782709059e-42\\
	0.225	-7.99419216079741e-34\\
	0.23125	-5.40341601271852e-29\\
	0.2375	-2.00915949017361e-19\\
	0.24375	-1.09552171449986e-12\\
	0.25	-1.09631660797738e-08\\
	0.25625	-2.45478937939807e-06\\
	0.2625	-4.56685834772883e-05\\
	0.26875	6.70385063680861e-05\\
	0.275	0.00157663760518493\\
	0.28125	0.00441138783235335\\
	0.2875	0.00720056157416548\\
	0.29375	0.00944195755063167\\
	0.3	0.0110343733636755\\
	0.30625	0.0120331921259723\\
	0.3125	0.0125715802250735\\
	0.31875	0.0127914914392327\\
	0.325	0.0128068367954391\\
	0.33125	0.0126979357755165\\
	0.3375	0.0125169599131595\\
	0.34375	0.0122968372814669\\
	0.35	0.0120577021511017\\
	0.35625	0.0118118631994467\\
	0.3625	0.0115666277006739\\
	0.36875	0.0113263284691771\\
	0.375	0.0110933764476887\\
	0.38125	0.0108690592126263\\
	0.3875	0.0106539211067243\\
	0.39375	0.0104480890232778\\
	0.4	0.0102514071737091\\
	0.40625	0.0100635772156128\\
	0.4125	0.00988420295422863\\
	0.41875	0.00971285397635297\\
	0.425	0.00954907718713074\\
	0.43125	0.00939242718520986\\
	0.4375	0.00924246592495335\\
	0.44375	0.00909877784317662\\
	0.45	0.00896096571646737\\
	0.45625	0.00882865843974434\\
	0.4625	0.00870150611079818\\
	0.46875	0.00857918410575105\\
	0.475	0.00846138846987358\\
	0.48125	0.00834783805488745\\
	0.4875	0.00823827053654378\\
	0.49375	0.00813244350733754\\
	0.5	0.00803013115428008\\
};
\addlegendentry{$\lambda_2$}

\end{axis}
\end{tikzpicture}
%
%%% lam 3 %%%
%
\begin{tikzpicture}

\begin{axis}[%
width=1.5in,
height=1.5in,
at={(0.758in,0.481in)},
scale only axis,
xmin=0,
xmax=0.5,
xlabel style={font=\color{white!15!black}},
xlabel={time t},
ymin=-0.01,
ymax=0.002,
axis background/.style={fill=white},
title style={font=\bfseries},
axis x line*=bottom,
axis y line*=left,
xmajorgrids,
ymajorgrids,
legend style={legend cell align=left, align=left, draw=white!15!black, at={(0.35,0.20)}}
]
\addplot [color=color0, draw=none, mark=*]
table[row sep=crcr]{%
	0.00625	0\\
	0.0125	0\\
	0.01875	0\\
	0.025	0\\
	0.03125	0\\
	0.0375	0\\
	0.04375	2.63717568000376e-317\\
	0.05	1.69162838525987e-302\\
	0.05625	7.97755183008671e-289\\
	0.0625	2.40721473382546e-276\\
	0.06875	-1.73274625391132e-263\\
	0.075	2.05465373714289e-251\\
	0.08125	4.59718536722041e-238\\
	0.0875	1.15726342928179e-226\\
	0.09375	2.62652964699959e-217\\
	0.1	-2.78061325861959e-205\\
	0.10625	-1.05445874222383e-195\\
	0.1125	5.69957932493213e-185\\
	0.11875	-8.99153085713257e-175\\
	0.125	1.31388931800404e-166\\
	0.13125	2.09841606861136e-156\\
	0.1375	-3.51695920405493e-146\\
	0.14375	-3.638504410253e-138\\
	0.15	-3.98784684875198e-129\\
	0.15625	4.872416598765e-120\\
	0.1625	-2.82937643349736e-113\\
	0.16875	5.60591676889781e-103\\
	0.175	-1.02358759820285e-95\\
	0.18125	6.76337762886208e-87\\
	0.1875	-1.72878458804612e-79\\
	0.19375	2.13765965541485e-71\\
	0.2	-5.76006409845596e-64\\
	0.20625	2.50506259586854e-56\\
	0.2125	-1.09539885044821e-48\\
	0.21875	1.09520780954592e-41\\
	0.225	-9.39990452010201e-34\\
	0.23125	-6.35384636288148e-29\\
	0.2375	-2.36245360061189e-19\\
	0.24375	-1.28816168681223e-12\\
	0.25	-1.289925581884e-08\\
	0.25625	-2.96963586504045e-06\\
	0.2625	-7.28842940050179e-05\\
	0.26875	-0.000398044576194121\\
	0.275	-0.000863130498142521\\
	0.28125	-0.0012370514139082\\
	0.2875	-0.00158797726811867\\
	0.29375	-0.00199587173720274\\
	0.3	-0.002463737096019\\
	0.30625	-0.00296910663461052\\
	0.3125	-0.00348207177566681\\
	0.31875	-0.00398055902480081\\
	0.325	-0.00445096975042605\\
	0.33125	-0.00488754796541642\\
	0.3375	-0.0052887679072386\\
	0.34375	-0.00565585068845143\\
	0.35	-0.00599098719579216\\
	0.35625	-0.00629693392973041\\
	0.3625	-0.0065763742085507\\
	0.36875	-0.00683192927463964\\
	0.375	-0.00706593231645012\\
	0.38125	-0.0072805253618204\\
	0.3875	-0.00747757278862872\\
	0.39375	-0.00765875300070556\\
	0.4	-0.00782552170351826\\
	0.40625	-0.00797918219254516\\
	0.4125	-0.00812086829687161\\
	0.41875	-0.00825159530118952\\
	0.425	-0.0083722515577062\\
	0.43125	-0.00848363485590425\\
	0.4375	-0.00858644819126327\\
	0.44375	-0.0086813257247274\\
	0.45	-0.00876883064427269\\
	0.45625	-0.00884947375817958\\
	0.4625	-0.00892371243162276\\
	0.46875	-0.00899196396449487\\
	0.475	-0.00905460508251971\\
	0.48125	-0.00911198160960414\\
	0.4875	-0.00916440824500313\\
	0.49375	-0.00921217559318897\\
	0.5	-0.00925555008942534\\
};
\addlegendentry{$\lambda_3$}

\end{axis}
\end{tikzpicture}	
	\caption{Values of the three Lagrange multipliers obtained by the cG scheme with $r = 1$ for the quasilinear heat equation example from \cref{subsection:quasilinear_heat_equation} with $c_1 = 3$ and $c_2 = 1$.}
	\label{fig:resistence:lambda}
\end{figure}

The boundary and transmission conditions between the two subdomains are the following.
At $z = 0$, we impose non-homogeneous Dirichlet conditions
\begin{subequations}
	\label{eq:quasilinear_heat_equations_constraints}
	\begin{align}
		\label{eq:quasilinear_heat_equations_Dirichlet_conditions}
		g_1(u,t)
		&
		\coloneqq
		u(0,t) - 1 
		.
		\intertext{At $z = 1$, we have the nonlinear transmission conditions}
		\label{eq:quasilinear_heat_equations_transmission_condition_1}
		g_2(u,t)
		&
		\coloneqq
		\frac{\partial}{\partial n^-} \paren[big][]{u(1^-,t)^{c_1}} + \alpha \paren[auto](){u(1^-,t) - u(1^+,t)}
		,
		\\
		\label{eq:quasilinear_heat_equations_transmission_condition_2}
		g_3(u,t)
		&
		\coloneqq
		\frac{\partial}{\partial n^+} \paren[big][]{u(1^+,t)^{c_2}} + \alpha \paren[auto](){u(1^+,t) - u(1^-,t)}
		.
	\end{align}
\end{subequations}
Here $\alpha = 10$ denotes the heat transfer coefficient between the two subdomains. 
At $z = 2$, we impose homogeneous Neumann boundary conditions, which are already included in~\verb!K!. 
Notice that the solution~$u$ of \eqref{eq:quasilinear_heat_equations} is discontinuous at $z = 1$ and $u(1^+,t)$ denotes the right limit while $u(1^-,t)$ is the left limit.
Moreover, $\frac{\partial}{\partial n^-} = \frac{\partial}{\partial z}$ denotes the normal derivative at $z = 1$ with respect to the \eqq{$-$}~domain~$(0,1)$ and $\frac{\partial}{\partial n^+} = -\frac{\partial}{\partial z}$ is the analog for the \eqq{$+$}~domain~$(1,2)$.
The conditions in \eqref{eq:quasilinear_heat_equations_constraints} are thus implemented as follows:
\begin{lstlisting}[frame=none,backgroundcolor=]
	g(1) = x(1) - 1;
	g(2) = (x(end/2+0).^c1 - x(end/2-1).^c1) / h;
	g(2) = g(2) + alpha * (x(end/2+0) - x(end/2+1));
	g(3) = (x(end/2+1).^c2 - x(end/2+2).^c2) / h;
	g(3) = g(3) + alpha * (x(end/2+1) - x(end/2+0));
\end{lstlisting}
The problem is solved on the time interval $[0,T]$ with $T = 0.5$ and time step size~$\Delta = T/40$. 
The initial condition is 
\begin{align*}
	u(x,0) 
	= 
	\begin{cases}
		1-4x& \text{in } (0,1/4), \\
		0& \text{else} 
	\end{cases}
\end{align*}
and hence consistent with the constraints, \ie, $g(u(z,0),0) = 0$.
Notice that both $f$ and $g$ are linear in case $c_1 = c_2 = 1$ and otherwise nonlinear.
An illustration showing some time steps of the solution using degree~$r = 1$ and exponents~$c_1 = 3$, $c_2 = 1$ is presented in \cref{fig:resistence:state}. 
Moreover, \cref{fig:resistence:lambda} shows the Lagrange multipliers.
One can clearly see the time when the temperature front reaches the coupling point between the two domains at $z = 1$, around time~$t = 0.25$.

For this example, an exact solution is not known.
Therefore, we estimate the discretization error using the numerical solution on a fine grid.
In the linear case, we can observe the convergence rates predicted by \cref{thm:convergenceDAE}.
In the nonlinear case with $c_1 = 3$ and $c_2 = 1$, we observe slightly reduced rates.

%
%=============================================================================
\subsection{Pendulum}\label{sec:numerics:pendulum}
In this final numerical experiment, we investigate whether the proposed cG scheme~\eqref{eq:cG_time_step} may also be applied to constrained Hamiltonian systems which have a similar structure, namely 
\begin{align*}
	J \dot x 
	& 
	= 
	-\nabla_x E(x,t) - g_x(x,t)^\transp \lambda
	, 
	\\
	0 
	& 
	= 
	g(x,t) 
	,
\end{align*}
with a skew-symmetric matrix~$J$. 
DAEs of this form typically appear in the field of mechanical systems and, in contrast to~\eqref{eq:DAE}, they have a differentiation index~3, cf.~\cite[Ch.~VII]{HairerWanner:1996:1}. 

As a prototypical example, we consider the first-order formulation of the mathematical pendulum with variables~$x = [x_1, x_2, y_1, y_2]^\transp$ and  
\begin{equation*}
	J 
	= 
	\begin{bmatrix} & & 1 & 0 \\  & & 0 & 1 \\ -1 & 0 & & \\ 0 & -1 & & \end{bmatrix}
	, 
	\qquad
	E(x) 
	= 
	\frac12\, y_1^2 + \frac12\, y_2^2 + \gamma\, x_2
	, 
	\qquad
	g(x) 
	= 
	x_1^2 + x_2^2 - \ell^2
	.
\end{equation*}
Here, $x_1, x_2$ denote the horizontal and vertical positions, $y_1, y_2$ are the respective velocities, $\ell$ denotes the length of the pendulum, and~$\gamma$ is the gravitational constant. 
Due to the matrix~$J$ in front of the differential term, the cG scheme~\eqref{eq:cG_time_step_simplified} needs to be adjusted accordingly. 
This simply means that the matrix~$J$ is added in front of the first sum in equation~\eqref{eq:cG_time_step_simplified_1}. 
Let us mention that a similar cG approach was considered in~\cite{EggerHabrichShashkov:2020:1} but applied to an index-reduced formulation. 
Yet another cG scheme for constrained mechanical systems (also after index reduction) was already introduced in~\cite{BetschSteinmann:2002:1}. 

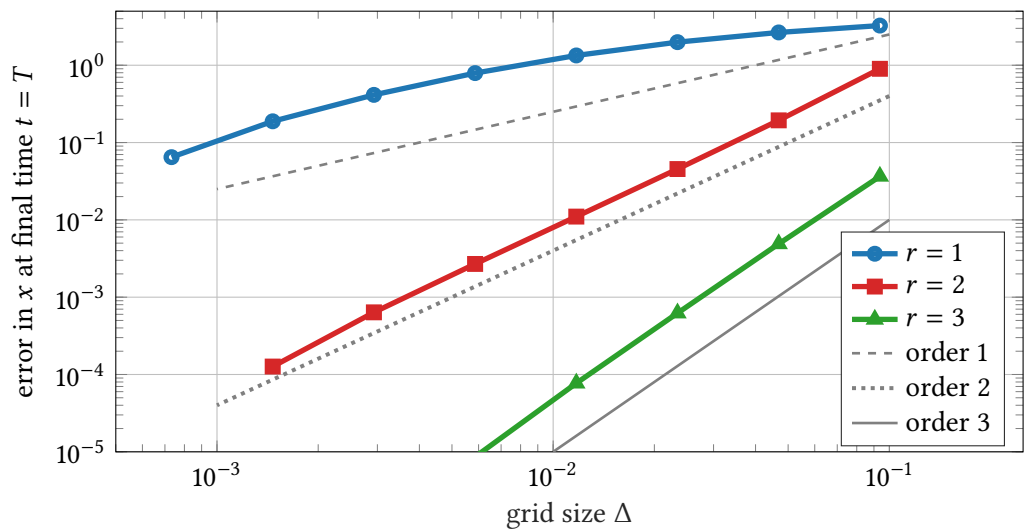
\begin{figure}
	\centering

% color for tikz images
\definecolor{color0}{rgb}{0.12156862745098,0.466666666666667,0.705882352941177}
\definecolor{color1}{rgb}{1,0.498039215686275,0.0549019607843137}
\definecolor{color2}{rgb}{0.172549019607843,0.627450980392157,0.172549019607843}
\definecolor{color3}{rgb}{0.83921568627451,0.152941176470588,0.156862745098039}
\definecolor{color4}{rgb}{0.580392156862745,0.403921568627451,0.741176470588235}

\begin{tikzpicture}

\begin{axis}[%
width=4.7in,
height=2.3in,
scale only axis,
xmode=log,
xmin=0.0005,
xmax=0.25,
xminorticks=true,
xlabel style={font=\color{white!15!black}},
xlabel={grid size $\Delta$},
ylabel={error in $x$ at final time~$t=T$},
ymode=log,
ymin=0.00001,
ymax=5,
yminorticks=true,
title style={font=\bfseries},
xmajorgrids,
%xminorgrids,
ymajorgrids,
%yminorgrids,
legend style={legend cell align=left, align=left, draw=white!15!black, at={(0.98,0.5)}}
]

\addplot [color=color0, line width=2.0pt, mark=o]
table[row sep=crcr]{%
%0.1875	3.66034632143634\\
0.09375	3.25883742700225\\
0.046875	2.65606527706962\\
0.0234375	1.9889645211368\\
0.01171875	1.33810459960557\\
0.005859375	0.791365703868973\\
0.0029296875	0.413716472058601\\
0.00146484375	0.188298829158363\\
0.000732421875	0.0647354268040043\\
};
\addlegendentry{$r=1$}

\addplot [color=color3, line width=2.0pt, mark=square*]
table[row sep=crcr]{%
%0.1875	4.78234980797966\\
0.09375	0.899127110200189\\
0.046875	0.193742567198494\\
0.0234375	0.0454352782101015\\
0.01171875	0.0110242624387597\\
0.005859375	0.00269053197129799\\
0.0029296875	0.000636666494910904\\
0.00146484375	0.000126694890820418\\
};
\addlegendentry{$r=2$}

\addplot [color=color2, line width=2.0pt, mark=triangle*]
table[row sep=crcr]{%
%0.375	8.13008969984945\\
%0.1875	0.34335419224785\\
0.09375	0.0365396265824078\\
0.046875	0.00489232574789361\\
0.0234375	0.000622117004094698\\
0.01171875	7.71329094778456e-05\\
0.005859375	8.6104650445551e-06\\
};
\addlegendentry{$r=3$}

%%%%%%%% order lines %%%%%%%%

\addplot [color=gray, line width=1pt, dashed]
table[row sep=crcr]{%
	0.1	2.5\\
	0.001	0.025\\
};
\addlegendentry{order $1$}

\addplot [color=gray, line width=1.5pt, dotted]
table[row sep=crcr]{%
	0.1	0.4\\
	0.01	0.004\\
	0.001	0.00004\\
};
\addlegendentry{order $2$}

\addplot [color=gray, line width=1pt]
table[row sep=crcr]{%
	0.1	0.01\\
	0.01	0.00001\\
};
\addlegendentry{order $3$}

\end{axis}
\end{tikzpicture}%
	\caption{Convergence orders for the state obtained by the cG scheme applied to the DAE formulation of the pendulum example from \cref{sec:numerics:pendulum}.}
	\label{fig:pendulum}
\end{figure}

The outcome of the numerical experiment is presented in \cref{fig:pendulum} and shows the convergence orders~$r$ of the state for this index-3 example with nonlinear constraint. 
Furthermore, the energy~$E$ is conserved up to an error of order~$r$ as well. 
We would like to emphasize, however, that the theoretical results of this paper do not apply to the present case. 
The reason is that in this setting $\kronD$ turns into~$\kronD = \overline D \otimes J$, which is still invertible, but the argumentation in the proof of~\Cref{proposition:invertibility_Newton_matrix} is no longer valid and we cannot deduce the invertibility of the matrix~$\kronZ^\transp \kronD \, \kronZ = \overline D \otimes (\overline Z^\transp J \overline Z) + \cO(\varepsilon)$.  
%
%
%=============================================================================
%=========  Conclusion
%=============================================================================
\section{Conclusion}
In this paper, we introduced and analyzed cG schemes for the numerical solution of semi-explicit DAEs of index~2. 
Based on the choice of the quadrature formula for the ODE part, described by~$f$, and the points where the constraint~$g$ is enforced, the assumed \emph{variational consistency} automatically determines the discretization of the Lagrange multiplier. 
For DAEs with linear constraints, we showed that~$x$ and~$\lambda$ converge to the exact solution with order~$r+1$ if the ansatz space consists of piecewise polynomials of order~$r$. % globally continuous,
Moreover, the numerical experiment in \cref{sec:numerics:pendulum} demonstrated the potential of the proposed cG schemes for the application to semi-explicit DAEs of index-3 as they appear in the modeling of mechanical systems. 

% \section*{TODO}
% \begin{itemize}
%   \item remove commented parts
%   \item make labels consistent
%   \item code repository
%   \item \LaTeX
% \end{itemize}

% Insert bibliography
\printbibliography

\end{document}